\documentclass[11pt]{article}
\usepackage{geometry} 
\usepackage{amsfonts}
\usepackage{epsfig}
\usepackage{url}
\usepackage{amssymb,latexsym,amsmath,amsthm,verbatim}
\usepackage{graphicx,epsfig,epstopdf,amssymb,color}
\usepackage{bm}
\newcommand{\N}{{\mathbb N}}

\newcommand{\veps}{\varepsilon}
\newcommand{\R}{\mathbb{R}} 
\newcommand{\C}{\mathbb{C}}
\newcommand{\K}{\mathbb{K}}
\newcommand{\OO}{\mathcal{O}}

\newcommand{\ba}{\begin{array}}
\newcommand{\ea}{\end{array}}

\newcommand{\eps}{\varepsilon}
\newcommand{\sign}{{\rm sign}}
\renewcommand{\epsilon}{\varepsilon}

\newtheorem{theorem}{Theorem}[section]
\newtheorem{lemma}[theorem]{Lemma}

\newtheorem{corollary}[theorem]{Corollary}
\newtheorem{remark}[theorem]{Remark}

\DeclareMathOperator*{\k0}{K_0 \!}
\DeclareMathOperator*{\kk}{K_0^{-1}\!}
\DeclareMathOperator*{\res}{Res }
\DeclareMathOperator*{\supp}{supp}
\DeclareMathOperator*{\sech}{sech}

\begin{document}
\title{Existence of long time solutions and validity of the Nonlinear Schr\"odinger approximation  
for a quasilinear dispersive equation}
\author{Wolf-Patrick D\"ull, Max He{\ss}}
\date{\today}
\maketitle

\begin{abstract}
We consider a nonlinear dispersive equation with a quasilinear quadratic term.
We establish two results. First, we show that solutions to this equation with initial data of order $\mathcal{O}(\varepsilon)$ in Sobolev norms exist for a time span of order $\mathcal{O}(\varepsilon^{-2})$ for sufficiently small $\varepsilon$. Secondly, we derive the  Nonlinear Schr\"odinger (NLS) equation as a formal approximation equation describing slow spatial and temporal modulations of the envelope of an underlying carrier wave, and justify this approximation with the help of error estimates in Sobolev norms between exact solutions of the quasilinear equation and the formal approximation obtained via the NLS equation.\\
The proofs of both results rely on estimates of appropriate energies whose constructions are inspired by the method of normal-form transforms. To justify the NLS approximation, we have to overcome additional difficulties caused by the occurrence of resonances.
We expect that the method developed 
in the present paper will also allow to prove the validity of the NLS approximation for a larger class of quasilinear dispersive systems with resonances.
\end{abstract}

\section{Introduction}
\label{sec1}
In this paper, we consider the quasilinear dispersive equation
\begin{equation}
\label{DGL}
 \partial_t u = \k0 u - u\, \partial_x u \,,
\end{equation}
where $ x,t, u(x,t) \in \R $ and
the linear operator $\k0\,$ is defined by its symbol 
\begin{equation}
\widehat{K}_0 \, (k)= - i \tanh(k)\,.
\end{equation} 

First, we show that solutions of \eqref{DGL} with initial data of order $\mathcal{O}(\varepsilon)$ in Sobolev norms exist for a time span of order $\mathcal{O}(\varepsilon^{-2})$ for sufficiently small $\varepsilon$, although equation \eqref{DGL} has a quadratic nonlinearity. More precisely, we prove 

\begin{theorem}
\label{Theorem 1 Langzeit loesung}
Let $s \ge 2$. There are constants $a, \varepsilon_0>0$ and $C\ge0$ such that
for all $\varepsilon \in\, (0, \varepsilon_0)$ and  $ u_0 \in H^s$ with
\begin{align*}
\Vert u_0 \Vert_{H^2} \le \varepsilon\,,
\end{align*}
there exists a solution $u \in C(I_\varepsilon, H^s)\, \cap\, C^1(I_\varepsilon, H^{s-1}) $, where $I_\varepsilon = [-a/\varepsilon^2, a/\varepsilon^2]$,
of \eqref{DGL} with $u(x,0) = u_0(x)$ for all $x \in \R$, which satisfies
\begin{align*}
\sup_{t \in I_\varepsilon}\Vert u(t) \Vert_{H^s} 
\le C \Vert u_0 \Vert_{H^s} \,.
\end{align*}
\end{theorem}

\begin{figure}[htbp]
\epsfig{file=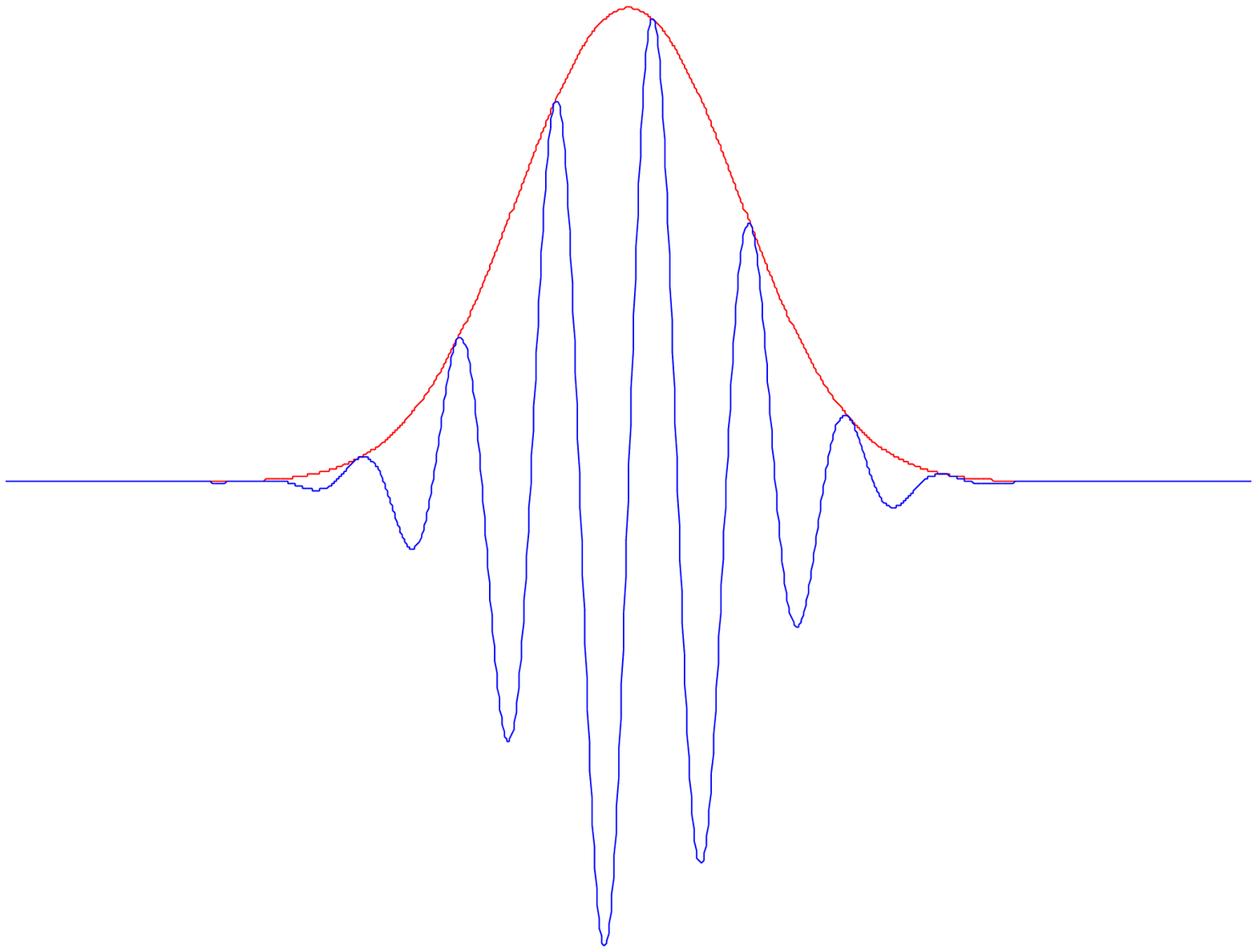,width=14.6cm,height=6.2cm,angle=0}
\vspace*{-6.3cm}

\hspace{1.9cm}
\epsfig{file=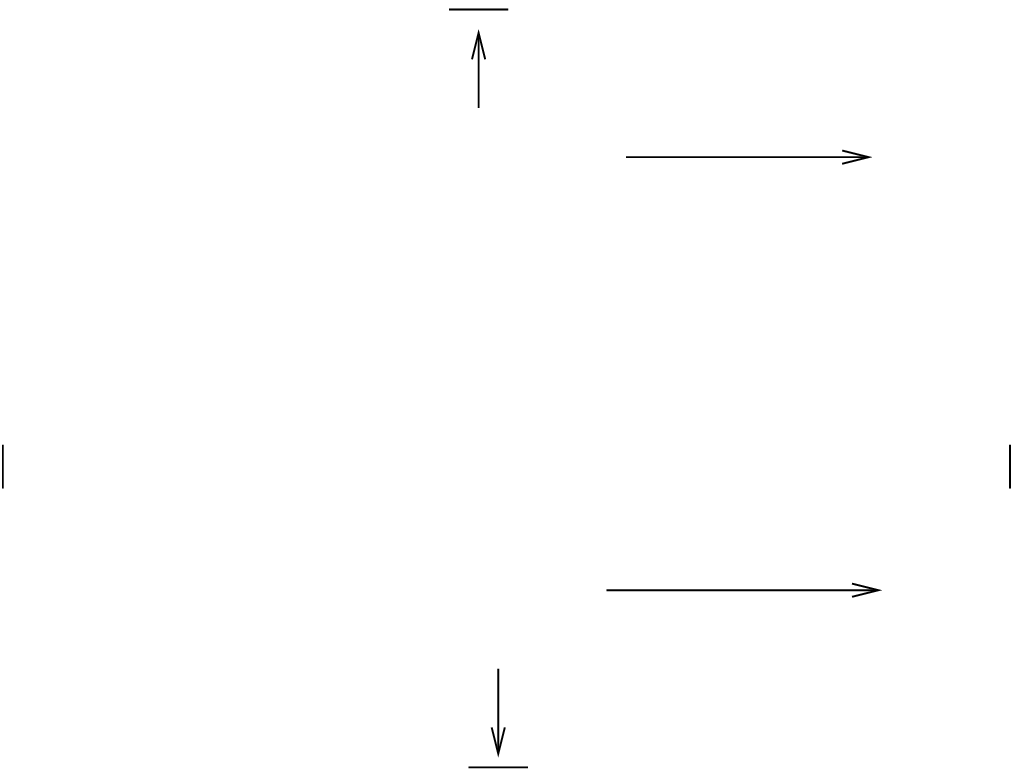,width=10.7cm,height=6.5cm,angle=0}

\vspace*{-5.34cm}
\hspace*{9.7cm}$c_\mathrm{g}$  \vspace*{3.2cm}

\hspace*{9.7cm} $c_\mathrm{p}$

\vspace*{-2.5cm}
\hspace{2.2cm}\hspace{5.0cm}$\varepsilon$
\vspace*{0.3cm}

\hspace{2.65cm}$1/\varepsilon$
\hspace{2.7cm}
\vspace*{2.6cm}

\caption{{\small The envelope (advancing with the group velocity 
$ c_g $) of the oscillating wave packet 
(advancing with the phase velocity 
$ c_p= \omega_0/k_0  $) is described by the 
amplitude $ A $ which solves  the NLS equation \eqref{NLS}.} \label{fig1} }
\end{figure}

Secondly, we derive the Nonlinear Schr\"odinger (NLS) approximation for equation \eqref{DGL} and prove its validity. The NLS equation plays an important role in describing  approximately  
slow modulations in time and space of an underlying spatially and temporarily oscillating wave packet in dispersive systems, for example, the water wave equations, see \cite{AS81}.
In order to derive the NLS approximation, we make the ansatz $u= \veps \psi = \veps \psi_{NLS} + \OO(\veps^2)$, with
\begin{equation} \label{NLS-ansatz}
\veps \psi_{NLS}(x,t) 
= \veps A(\eps
(x-c_g t),\eps^2t) e^{i( k_0 x - \omega_0 t)} + \mathrm{c.c.} \,.
\end{equation}
Here $0 < \eps \ll 1$ is a small perturbation parameter,
$ \omega_0 > 0$ the basic temporal 
wave number associated to the basic spatial wave number $ k_0 > 0$ of the underlying carrier wave $ e^{i(k_0 x - \omega_0 t)}$,
$c_g$ the group velocity, $A$ the complex-valued amplitude, and c.c. the complex conjugate. 
With the help of (\ref{NLS-ansatz}) we describe slow spatial and temporal modulations 
of the envelope of the underlying carrier wave.
Inserting the above ansatz into \eqref{DGL} we find that $A$ satisfies at leading order in $\eps$ the NLS equation
\begin{equation} \label{NLS}
\partial_T A  = i \nu_1 \partial_{X}^2 A + i  \nu_2 A|A|^2\,, 
\end{equation}
where $X = \eps(x-c_g t)$, $ T = \eps^2 t $, and    $\nu_j= \nu_j(k_0) \in \R$. $ T $  is the slow time scale and $ X$
is the slow spatial scale, that means,
the time scale of the modulations is  
$\OO({1/\eps^2})$ and the spatial scale of the modulations
is $\OO({1/\eps})$. See Figure \ref{fig1}.
The basic spatial wave number $k= k_0$ and the basic temporal
wave number $\omega = \omega_0$ 
are related via  the linear dispersion relation
of equation (\ref{DGL}), namely 
\begin{equation}
\label{lindis}
\omega(k) = \tanh(k)\,.
\end{equation}
Then the group velocity $ c_g $ of the wave packet
is given by $ c_g = \partial_k  \omega|_{k=k_0}$.
Our ansatz leads to waves moving to the right. To obtain waves moving
to the left, $-\omega_0$ and $c_g$ have to be replaced by  $\omega_0$ and $-c_g$.

\vspace*{0.35cm}

To justify the NLS approximation for \eqref{DGL}, we prove
\begin{theorem}
\label{Theorem 2 NLS-APP}
Fix $s_A \geq 7$.  Then for all $k_0 > 0$ 
and for all $C_1,T_0 > 0$ there exist $C_2 > 0$, $\eps_0 > 0$ 
such that for all solutions $A \in
C([0,T_0],H^{s_A}(\R,\C))$ of the NLS equation (\ref{NLS})
with 
$$ 
\sup_{T \in [0,T_0]} \| A(\cdot,T) \|_{H^{s_A}(\R,\C)} \leq C_1
$$  
the following holds.
For all $\eps \in\, (0,\eps_0)$
there are solutions 
$$
u \in
C([0,T_0/\eps^2], H^{s_A}(\R,\R))
$$
of equation \eqref{DGL} which satisfy
  $$\sup_{t\in[0,T_0/\eps^2]} \| u(\cdot,t) -
    \eps \psi_{NLS}(\cdot,t)\|_{H^{s_A}(\R,\R)} 
\le C_2 \eps^{3/2}.$$
\end{theorem}

The error of order $\OO(\eps^{3/2})$ is small compared with the solution 
$u$ and the approximation $\eps \psi_{NLS}$, which are both of order
$\OO(\eps)$ in $ L^{\infty} $ such that the dynamics 
of the NLS equation can be found in equation \eqref{DGL}, too. The NLS equation is a completely integrable
Hamiltonian system, which can be solved 
explicitly with the help of some inverse scattering scheme, see, for example, \cite{AS81}.

It should be noted that the smoothness in our error bound is equal to the assumed smoothness of the amplitude. This can be achieved by using a modified approximation which has compact support in Fourier space but differs only slightly from $\eps \psi_{NLS}$. Such an approximation can be constructed because the Fourier transform of $\eps \psi_{NLS}$ is sufficiently strongly concentrated around the wave numbers $\pm k_0$.

We remark that such an approximation theorem should not be taken for granted. There are various counterexamples,
where approximation equations derived by reasonable formal arguments
make wrong predictions about the dynamics of the 
original systems, see, for example, \cite{Schn05,SSZ15}. For an introduction into theory and applications of the NLS approximation we refer to \cite{Schn11OWbuch}.
\medskip

Now, we explain the main ideas for the proofs of our theorems. Like in many other proofs of related estimates in the literature we will assume in our proofs of Theorem \ref{Theorem 1 Langzeit loesung} and Theorem \ref{Theorem 2 NLS-APP} that $s$ and $s_A$ are integers in order to simplify the analysis by using Leibniz's rule,  but our proofs can be generalized to be valid for all $s \geq 2$ and $s_A \geq 6$.

The main difficulty in the proof of Theorem \ref{Theorem 1 Langzeit loesung} is
to show that $I_\varepsilon$ is of order $\mathcal{O}(\varepsilon^{-2})$. If $u_0$ is of order $\mathcal{O}(\varepsilon)$, then, due to the fact that the nonlinear term $-u \partial_x u$ is quadratic, direct energy estimates only guarantee an existence interval of order $\mathcal{O}(\varepsilon^{-1})$ for 
$u$. A standard strategy to address this problem is to try to eliminate the quadratic term and transfer it into a cubic term with the help of a normal-form transform of the form
\begin{equation}
\label{nft}
\tilde{u}:= u + N(u,u)\,,
\end{equation}
where $N$ is an appropriately constructed bilinear mapping, see \cite{Sh85, Kal88}. In the case of equation \eqref{DGL}, a direct computation of the evolution equation for $\tilde{u}$ with the help of equation \eqref{DGL} yields
that $\tilde{u}$ solves 
an evolution equation of the form
\begin{equation} 
\label{t-ueq}
\partial_t \tilde{u} = \k0 \tilde{u} + h(u, \partial_x u)\,,
\end{equation}
where $h(u, \partial_x u)$ is a cubic term if $N$ satisfies
\begin{equation} 
\label{Nrel}
-\k0 N(u,u) + N(\k0 u, u) + N(u, \k0 u) = \frac{1}{2} \partial_x (u^{2})\,.
\end{equation}
Since $\k0\,$ satisfies the identity
\begin{align}
\label{K0id1}
\k0 \, (fg)-\k0 \,(f)\,g-f\k0 \,(g)=\k0 \;(\k0 \, (f) \k0 \,(g))\,,
\end{align}
see Lemma \ref{k0-id-lem} below, it follows
\begin{equation} 
\label{Ndef1}
N(u,u) = -\frac{1}{2} \kk \partial_x (\kk u)^{2} \,.
\end{equation}
However, this condition for $N$ causes two problems. The first problem is that $\kk u$ may not exist, and the second one is that $N(u,u)$ loses one derivative, that means, $ u \mapsto N(u,u)$ maps $H^{m+1}(\R,\C)$ into $H^{m}(\R,\C)$ or $C^{n+1}(\R,\C)$ into $C^{n}(\R,\C)$. Even if it was possible to invert the normal-form transform \eqref{nft}, the cubic term $h$ expressed in terms of $\tilde{u}$ would lose two derivatives such that it would not be possible to use equation \eqref{t-ueq} to derive closed energy estimates for $\tilde{u}$.

To overcome these problems, we do not perform the normal-form transform \eqref{nft} explicitly, but only use the term $N$ to construct an energy of 
the form
\begin{align}
\label{Edefa}
\mathcal{E}_{s}= \sum_{\ell=0}^{s} E_{\ell}\,,
\end{align}
where the summands $E_{\ell}$ are defined by a slight, $\ell$-dependent modification of the equation 
\begin{align}
\label{Edefb}
E_{\ell} = \frac{1}{2} \Vert \partial_x^{\ell} u \Vert_{L^2}^2 +  \int_{\R}     
\partial_x^{\ell} u\, \partial_x^{\ell}N(u,u)\,dx
\end{align}
to get around the problem that $\kk u$ may not exist. More precisely, since
\begin{equation} 
\frac{d}{dt} \|u\|_{L^{2}} = 0
\end{equation}
for any sufficiently regular solution $u$ of \eqref{DGL}, see Lemma \ref{E0-lem} below, we define
\begin{align}
\label{E0-def}
E_{0}:= \frac{1}{2} \Vert u \Vert_{L^2}^2\,.
\end{align}
Moreover, due to
\begin{align}
\label{partintf}
\int_{\R}     
\partial_x^{\ell} f\, \partial_x^{\ell}(f\,\partial_x f)\,dx
&= \sum_{a=1}^{\ell-1} {\binom {\ell} a} 
\int_{\mathbb{R}}
\partial_x^{\ell} f
\, \partial_x^a f\, \partial_x^{\ell-a+1} f\, dx 
+
\frac{1}{2}
\int_{\mathbb{R}}
\partial_x^{\ell} f\,
\partial_x^{\ell} f\, \partial_x f\,dx
\end{align}
for sufficiently regular functions $f$, which follows with the help of Leibniz' rule and integration by parts, and because of the facts that $\kk$ is skew symmetric and $\kk \partial_x g$ exists for any $g \in H^{1}$, we define
\begin{align}
\label{Eell-def}
E_{\ell} :=& \;\frac{1}{2} \| \partial_x^{\ell} u\|_{L^2}^2  
+\sum_{a=1}^{\ell-1} {\binom \ell a} \int_{\mathbb{R}} \kk \partial_x^{\ell} u \kk \partial_x^{a} u \kk \partial_x^{\ell-a+1} u\, dx
\\[2mm]
& \nonumber \;
+ \frac{1}{2} \int_{\mathbb{R}} \kk \partial_x^{\ell}u \kk \partial_x^{\ell} u \kk \partial_x u \,dx
\end{align}
for $\ell >0$. 

$\sqrt{\mathcal{E}_{s}}$ is equivalent to $\Vert u \Vert_{H^{s}}$ for $s \ge 2$ and $\Vert u \Vert_{H^{1}} = \mathcal{O}(\varepsilon)$, see Lemma \ref{lem Energie aequi}.
Due to the skew symmetry of $\k0$ and \eqref{Nrel}, the right-hand side of the evolution equation for $\mathcal{E}_{s}$ contains neither quadratic nor cubic terms. Moreover, the right-hand side of the evolution equation for $\mathcal{E}_{s}$ can be written as a sum of integral terms containing at most one factor $\partial_x^{s+1} u$ and not two. 
Consequently, using integration by parts and estimates for the commutator 
$[K_0^{-1}, u] \partial_x u$, we obtain
\begin{equation}
\frac{d}{dt} \mathcal{E}_s \lesssim \veps^{2} {\mathcal{E}_s}
\end{equation}
as long as $\|u\|_{H^2} = \OO(\varepsilon)$ such that Gronwall's inequality yields the $\mathcal{O}(1)$-boundedness of ${\mathcal{E}_s}$ and hence of $u$ for all $t\in I_\varepsilon$. For further details, see Section \ref{sec2}.
\medskip

There is an equation which is related to \eqref{DGL}, namely
\begin{equation}
\label{DGL-H}
 \partial_t u = {\rm H} u - u\, \partial_x u \,,
\end{equation}
where $ x,t, u(x,t) \in \R $ and ${\rm H}$ is the Hilbert transform. For this equation, the analog of 
Theorem \ref{Theorem 1 Langzeit loesung} was proven in \cite{HITW13}. The proof also relies on energy estimates inspired by a normal-form transform of the form \eqref{nft}, but the details of the proof are simpler in the following sense. Since the Hilbert transform also satisfies the identity \eqref{K0id1} (with $\k0\,$  replaced by ${\rm H}$), one obtains for the bilinear mapping $N$ the condition \eqref{Ndef1} with $\kk\,$
replaced by ${\rm H}^{-1}$. Because ${\rm H}^{-1}=-{\rm H}$ is well-defined in $L^2$, an appropriate energy can be defined directly by \eqref{Edefa} and \eqref{Edefb}.

In \cite{HIT14}, \cite{IT14} and \cite{IT16}, the techniques from \cite{HITW13} were further developed and applied to the 2D water wave problem with infinite depth in holomorphic coordinates in order to derive
high-order energy estimates which allowed the authors  
to prove local well-posedness in Sobolev spaces and to establish extended life spans for small solutions. Moreover, by combining those high-order energy estimates with dispersive decay estimates the authors showed global existence of small localized solutions.
\medskip

In order to prove Theorem \ref{Theorem 2 NLS-APP}, we estimate the error
\begin{align}
\label{R1}
\eps^{\beta} R := u- \eps \psi
\end{align}
 for all $t \in [0, T_0/\eps^2]$ to be of order $\mathcal{O}(\eps^{\beta})$ in $H^{s_A}(\R,\R)$ for a $\beta \geq 3/2$, that means, we prove that
$R$ is of order $\OO(1)$ for all $t \in [0,T_0/\eps^2]$. The error $R$ satisfies the equation
\begin{align}
\label{erreq1}
\partial_t R= \k0 R +2 \eps B(\psi, R) + \eps^{\beta} B(R, R) +\varepsilon^{-\beta} \res(\eps \psi)
\end{align}
with
\begin{align}
\label{B}
B(f, g) = -\frac{1}{2}\partial_x(fg)\,.
\end{align}
If $A$ is a sufficiently regular solution of the NLS equation \eqref{NLS}, then the Fourier transform $\veps \widehat{\psi}_{NLS}$ is so strongly concentrated around the wave numbers $\pm k_0$ that it is possible to construct an approximation function $ \eps \psi$ with compact support in Fourier space satisfying $\eps \psi = \veps \psi_{NLS} + \OO(\veps^{3/2})$ in $H^{s_A}(\R,\R)$ and to choose
$\beta$ such that 
\begin{align}
\partial_t R= \k0 R +2 \eps B(\psi,R) + \mathcal{O}(\varepsilon^{2})
\end{align}
with respect to a Sobolev norm.
Moreover, the approximation $ \eps \psi$ can be split into 
\begin{align}
\eps \psi = \eps \psi_c + \eps^{2} \psi_s
\end{align}
with
\begin{align}
\label{supppsic}
\supp \widehat{\psi}_{c} = \{k \in \R:\, |k \mp k_0| \leq \delta\}\,,
\end{align}
where $\delta \in\, (0, k_0)$ is small, but independent of $\veps$, and $\psi_s = 
\mathcal{O}(1)$ with respect to a suitable norm, see Lemma \ref{Res-abschaetzung} below. Therefore, we have
\begin{align}
\partial_t R= \k0 R +2 \eps B(\psi_c,R) + \mathcal{O}(\varepsilon^{2})
\end{align}
such that the main difficulty is to control the quadratic term $2\eps B({\psi}_c,R)$ for 
$t \in [0, T_0/\eps^2]$.

It is again instructive to try to eliminate this term
with the help of a normal-form transform. 
We note that because $\psi_s$ is of order $\mathcal{O}(1)$, the other component $2 \eps^2 B(\psi_s,R)$ of $2 \eps B(\psi,R)$ is only of order $\mathcal{O}(\varepsilon^{2})$ and need not to be eliminated, which will simplify the construction of the normal-form transform significantly.

Hence, we look for a normal-form transform of the form
\begin{equation} \label{inft} 
\tilde{R} := R + \eps
N({\psi}_c,R)\,,
\end{equation}
with an appropriate bilinear mapping $N$, to obtain
\begin{align}
\label{tReq}
\partial_t \tilde{R}= \k0 \tilde{R} + \mathcal{O}(\varepsilon^{2})\,,
\end{align}
which yields
\begin{equation} 
\label{Nrel2}
-\k0 N({\psi}_c,R) + N(\k0 {\psi}_c, R) + N({\psi}_c, \k0 R) = -2B({\psi}_c, R)
\end{equation}
such that because of \eqref{K0id1} and \eqref{B} it follows
\begin{equation} 
\label{NN}
N({\psi}_c,R) = -\kk \partial_x (\kk \psi_c \kk R) \,.
\end{equation}
$\kk \psi_c$ exists due to \eqref{supppsic}, but
we have the problems that $\kk R$ may not exist and 
that $R \mapsto N({\psi}_c,R)$ loses one derivative.  
Since the $L^{2}$-norm of $R$ is not a conserved quantity and $N$ depends on the two different functions ${\psi}_c$ and $R$, we cannot define an energy in an analogous way as in the proof of Theorem \ref{Theorem 1 Langzeit loesung} to overcome these problems. Nevertheless, it is still possible to use the method of normal-form transforms for constructing an appropriate energy to control the error, but it takes some additional effort.

The problem that $\kk R$ may not exist is related to the occurrence of so-called resonances. In Fourier space, we have 
\begin{align}
\widehat{N}(\psi_c,R)(k)= \int_{\R}\widehat{n}(k, k-m,m)\, \widehat{\psi}_c(k-m)\,\widehat{R}(m)\, dm
\end{align}
with
\begin{align}
\widehat{n}(k,k-m,m)=- ik\,\widehat{K}_0^{-1}(k)\,\widehat{K}_0^{-1}(k-m)\, \widehat{K}_0^{-1}(m)\,.
\end{align}
Because of \eqref{supppsic}, it is instructive to analyze the behavior of $\widehat{n}$ for $|k-m| \approx k_0$. We have
\begin{align}
\label{nfrac}
\widehat{n}(k,k-m,m) \approx  \left\{ \begin{array}{ll}
 - \dfrac{ik}{\widehat{K}_0(k)\,\widehat{K}_0(k_0)\, \widehat{K}_0(k - k_0)}
 & \quad {\rm for }\;\, k-m \approx k_0\,, \\[5mm]
 - \dfrac{ik}{\widehat{K}_0(k)\,\widehat{K}_0(- k_0)\, \widehat{K}_0(k + k_0)}
 & \quad {\rm for }\;\, k-m \approx -k_0\,.
 \end{array} \right.
\end{align}
The denominators of the fractions in \eqref{nfrac} have the following zeros, which are called resonances.
Both denominators have a zero at $k=0$. Since the numerators also vanish at $k=0$ and $\lim_{|k|\to 0} k/\tanh(k) =1$, the singularity at $k=0$ is removable.
Such a resonance is called a trivial resonance. The fact that the resonance at
$k=0$ is trivial correlates with the fact that $\kk \partial_x g$ exists for any $g \in H^{1}$. Moreover, both denominators have one more zero - the first
denominator at $k=k_0$, the second one at $k=-k_0$. 
At these resonances, the respective numerators do not vanish. Such a resonance is called a non-trivial resonance. The fact that the resonances at
$k=\pm k_0$ are non-trivial correlates with the fact that $\kk R$ may not exist.

In the situation of a trivial resonance at $k=0$ and non-trivial resonances at 
$k=\pm k_0$, it is possible to apply a technique from \cite{DS06} for constructing a modified normal-form transform. The essential tools for the construction procedure from \cite{DS06} are as follows.

Since $\widehat{B}(\psi_c,R)(k)$ vanishes at $k=0$, one can expect that
$\widehat{R}(k)$ will grow for $k$ near $0$ more slowly than for $k$ further away from $0$. Hence, it makes sense to rescale the error with the help of the weight function
\begin{equation}
\widehat{\vartheta}(k) = \left\{ \begin{array}{cc} 1  &\quad {\rm for }\; |k| > \delta\, ,\\[1mm]
\epsilon + (1-\epsilon) |k|/\delta  &\quad {\rm for }\; |k| \leq \delta \, ,
\end{array} \right.
\end{equation}
where $\delta$ is chosen as above. More precisely, by writing
\begin{align}
\label{Rsc}
u= \eps \psi_c + \eps^{2} \psi_s + \eps^{5/2} \vartheta R\,,
\end{align}
where  $\psi_c$ and $\psi_s$ are as above and $\vartheta R$ is defined by $\widehat{\vartheta} \widehat{R}$, one obtains for the rescaled error $R$ an evolution equation of the form
\begin{equation}
\label{d/dt Rsc}
\partial_t R= \k0R - \eps \vartheta^{-1} \partial_x(\psi_c\, \vartheta {P}_{\varepsilon,\infty} R) + \mathcal{O}(\varepsilon^{2})\,.
\end{equation}  
Here, ${P}_{\varepsilon,\infty}$ is a linear operator with the symbol $\widehat{P}_{\varepsilon,\infty}(k)= (1-
\chi_{[-\eps,\eps]})(k)$, where $\chi_{[-\eps,\eps]}$ is the characteristic function on $[-\eps,\eps]$.
Now, constructing a normal-form transform of the form \eqref{inft} yields
\begin{equation}
\label{Nmod}
N(\psi_c, R) = - \vartheta^{-1} \kk \partial_x( \kk \psi_c\, \kk \vartheta P_{\varepsilon, \infty} R)\,,
\end{equation}
where $\kk \vartheta P_{\varepsilon, \infty} R$ exists for any $R \in L^{2}$.
However, since $(\widehat{\vartheta}(k))^{-1} = \mathcal{O}(\varepsilon^{-1})$
for $|k| < \delta$, the transformed 
error $\tilde{R}$ satisfies an evolution equation of the form
\begin{align}
\label{d/dt Rsct}
\partial_t \tilde{R}= \k0 \tilde{R} - \eps \sum_{j =\pm 1} (1- {P}_{\delta,\infty}) N(\psi_j, 2\eps \vartheta^{-1} \partial_x(\psi_j\, \vartheta {P}_{\varepsilon,\infty} R) ) + \mathcal{O}(\varepsilon^{2})\,,
\end{align}
with $\widehat{\psi}_j=\widehat{\psi}_c|_{[jk_0-\delta, jk_0+\delta]}$ and
$(1- {P}_{\delta,\infty}) N(\psi_j, 2\eps \vartheta^{-1} \partial_x(\psi_j\, \vartheta {P}_{\varepsilon,\infty} R))= \mathcal{O}(1)$.
But the term of order $\mathcal{O}(\varepsilon)$ on the right-hand side of \eqref{d/dt Rsct} can be eliminated with the help of a second normal-form transform of the form
\begin{align}
\label{nft2}
\check{R}= \tilde{R}+ \eps^{2} \sum_{j =\pm 1} \mathcal{T}_j(\psi_j,\psi_j,R)\,
\end{align}
with appropriate trilinear mappings $\mathcal{T}_j$. The construction of the trilinear mappings is similar to the construction of bilinear mappings for normal-form transforms. In the case of equation \eqref{d/dt Rsct}, no resonances occur
in the context of the construction of the trilinear mappings such that straightforward calculations yield 
\begin{equation}
\widehat{\mathcal{T}}_j(\psi_j,\psi_j,R)(k)= 
\int_{\mathbb{R}}\int_{\mathbb{R}}
\widehat{t}_j(k)\, \widehat{\psi}_j(k-m)\, \widehat{\psi}_j(m-n)\,\hat{R}(n)\, dn\, dm
\end{equation}
with 
\begin{align}
\widehat{t}_j(k)= &-\frac{k(k-jk_0)\,\widehat{\vartheta}(k-2jk_0)\, \chi_{[-\delta,\delta]}(k)}{\widehat{\vartheta}(k) 
\tanh(k) \tanh(jk_0) \tanh(k-jk_0)} \\[2mm] 
\nonumber & \times (\tanh(k) -2 \tanh(jk_0)-\tanh(k-2jk_0))^{-1} \,.
\end{align}
After these two normal-form transforms we have
\begin{align}
\label{d/dt Rcheck}
\partial_t \check{R}= \k0 \check{R} + \mathcal{O}(\varepsilon^{2})\,.
\end{align}
For further details about the two normal-form transforms discussed just now, we refer to \cite{DS06}.

However, since the error equation \eqref{erreq1} is quasilinear, also the modified normal-form transform $R \mapsto \check{R}(R)$ loses one derivative. It can be shown that this normal-form transform is nevertheless invertible, but the
term of order $\mathcal{O}(\varepsilon^{2})$ in the transformed error equation 
\eqref{d/dt Rcheck} loses two derivatives if it is expressed in terms of $\check{R}$.

To overcome the regularity problems, we pursue again the strategy from the proof of Theorem \ref{Theorem 1 Langzeit loesung} that we do not perform the normal-form transform explicitly, but only use it to construct an energy of 
the form
\begin{align}
\label{Edefa2}
\tilde{\mathcal{E}}_{s}= \sum_{\ell=0}^{s} \tilde{E}_{\ell}\,,
\end{align}
where the summands $\tilde{E}_{\ell}$ are defined by a slight, $\ell$-dependent modification of the equation 
\begin{align}
\label{Edefb2}
\tilde{E}_{\ell} = \frac{1}{2} \Vert \partial_x^{\ell} R \Vert_{L^2}^2 +  \epsilon \int_{\R}     
\partial_x^{\ell} R\, \partial_x^{\ell}N(\psi_c,R)\,dx +
\epsilon^{2} \sum_{j =\pm 1} \int_{\R}     
\partial_x^{\ell} R\, \partial_x^{\ell}\mathcal{T}_j(\psi_j,\psi_j,R)\,dx\,,
\end{align}
where $R$ is defined by \eqref{Rsc}, $N$ is defined by \eqref{Nmod}, and the 
mappings $\mathcal{T}_j$ are as in \eqref{nft2}.
Since $k^{\ell} (\widehat{\vartheta}(k))^{-1} = \mathcal{O}(1)$
for $|k| < \delta$ if $\ell \geq 1$,
we do not need to include the second normal-form transform in our energy for $\ell \ge 1$.  Hence, we define
\begin{align}
\tilde{E}_{\ell} := \frac{1}{2} \Vert \partial_x^{\ell} R \Vert_{L^2}^2 +  \epsilon \int_{\R}     
\partial_x^{\ell} R\, \partial_x^{\ell}N(\psi_c,R)\,dx
\end{align}
for $\ell \geq 1$. Then integration by parts yields
\begin{align}
\tilde{E}_{\ell} = \frac{1}{2} \Vert \partial_x^{\ell} R \Vert_{L^2}^2 +   \epsilon\, \mathcal{O}(\Vert R \Vert_{H^{\ell}}^2) 
\end{align}
for $\ell \geq 1$. Because the mapping
\begin{align*}
R \mapsto \int_{\mathbb{R}}  R \, \check{R}(R)\, dx 
\end{align*}
is in general not positive definite, we have to perform the full normal-form transform in the case of $\ell=0$ and define
\begin{align}
\tilde{E}_0 := \Vert \check{R}\Vert_{L^2}^2  \,.
\end{align}
The resulting loss of regularity does not mind here because it can be compensated 
with the help of the other components of our energy such that we obtain the equivalence of $\sqrt{\tilde{\mathcal{E}}_{s}}$ and $\Vert R \Vert_{H^s}$ for $s \geq 1$ and sufficiently small $\eps$, see Corollary \ref{en-Aequi}. Consequently, the right-hand side of the evolution equation of $\tilde{\mathcal{E}}_{s}$ can be written as a sum of integral terms containing at most one factor $\partial_x^{s+1} R$ and not two. Moreover, since $\|\check{R}\|_{H^s}^2$ differs from  $\tilde{\mathcal{E}}_{s}$ only by terms of order $\mathcal{O}(\veps^2)$, the evolution equations of $\tilde{\mathcal{E}}_{s}$ and $\|\check{R}\|_{H^s}^2$ share the property that their right-hand sides are of order $\mathcal{O}(\veps^2)$. Therefore, by using integration by parts, we obtain
\begin{equation}
\partial_t \tilde{\mathcal{E}_s} \lesssim \veps^{2} (\tilde{\mathcal{E}_s}+1)
\end{equation}
as long as $\|R\|_{H^s} = \OO(1)$ such that Gronwall's inequality yields the $\mathcal{O}(1)$-boundedness of $\tilde{\mathcal{E}}_s$ and hence of $R$ for all $t\in[0,T_0/\veps^{2}]$. 
\medskip

For the reasons discussed above, the justification of the NLS approximation for dispersive
systems with quasilinear quadratic terms is a highly nontrivial problem, which
has been remained unsolved in general for more than four decades. The first and very general NLS approximation theorem for quasilinear dispersive
wave systems was shown in \cite{Kal88}. However, the occurrence of quasilinear quadratic 
terms was excluded 
explicitly.
In the case of quasilinear quadratic terms, an NLS approximation theorem was proven for dispersive wave systems where the right-hand sides lose only half a derivative. The 2D water wave problem without surface tension and finite depth in Lagrangian coordinates falls into this class.
 In this case
the elimination of the quadratic terms is possible with the help of normal-form transforms. The right-hand sides of 
the transformed systems then lose one derivative and can be handled with the help of the Cauchy-Kowalevskaya theorem \cite{SW10,DSW12}.
Furthermore, the NLS approximation was justified for the 2D and 3D water wave problem without surface tension and infinite depth \cite{TW11,T14} by finding a different 
transform adapted to the special structure of that problem. 
Similarly, for 
the quasilinear Korteweg-de Vries
equation the result can be obtained by simply applying a
Miura transform \cite{Schn11}.
In \cite{CDS15}, the NLS approximation of time oscillatory long waves for equations with quasilinear quadratic terms was proven for analytic data without using a normal-form transform.
Moreover, another approach to address the problem of the validity of the NLS approximation can be found in \cite{MN13}.
Finally, some numerical evidence that the NLS approximation is also valid for quasilinear equations was given in \cite{CS11}.

Very recently, the first validity proof of the NLS approximation 
of a nonlinear Klein-Gordon equation with a quasilinear quadratic term in Sobolev spaces was given in \cite{D16}. The proof also relies on estimates of an appropriate energy which is constructed with the help of a normal-form transform. The construction of the energy is easier in the sense that no problems with resonances occur, but more difficult in the sense that the energy has to allow to control a system of two coupled error equations. 

Theorem \ref{Theorem 2 NLS-APP} of the present paper
is the first validity result for the NLS approximation of a quasilinear dispersive equation with resonances in Sobolev spaces.
\medskip

The plan of the paper is as follows.  In Section \ref{sec2} we prove Theorem \ref{Theorem 1 Langzeit loesung}. In Section \ref{sec3} we derive the NLS approximation. In Section \ref{sec4} we perform
the error estimates to prove Theorem \ref{Theorem 2 NLS-APP}. 
\medskip

The relevance of studying equation \eqref{DGL} lies in the fact that this equation serves as a simple model equation incorporating principal difficulties which have to be overcome both for 
establishing extended life spans for small solutions and for justifying the NLS approximation
for more complicated dispersive systems with resonances and rough nonlinearities. In particular, 
equation \eqref{DGL} and the 2D water wave problem with finite depth in various coordinates share the difficulties of having linear dispersion relations which cause a trivial resonance at the wave number $k=0$ as well as non-trivial resonances at $k= \pm k_0$ and possessing quadratic transport terms   
which preclude the application of the standard method of normal-form transforms because of a loss of regularity problem.

In the proofs of the results of the present paper,
we have further developed our approach from \cite{D16}, i.e., the replacement of the standard method of normal-form transforms by the use of an energy which includes essential parts of a normal-form transform, in the following sense. 
We have refined the basic form \eqref{Edefa}-\eqref{Edefb} of such an energy in a way that all problems caused by the occurring resonances can be
circumvented, which has taken some extra effort in the case of the justification of the NLS approximation. 

We expect that a combination of this type of energy with the types of energies we have constructed in \cite{D16, DS06} will yield the main component of an energy which will allow us to prove new justification
results for the NLS approximation of the 2D water wave problem with finite depth and other complicated dispersive systems with resonances.
In particular, we intend in forthcoming papers 
to establish an extended time span of the validity of the NLS approximation of the 2D water wave problem with finite depth and without surface tension as well as to solve the open problem of justifying the NLS approximation of the 2D water wave problem with finite depth and with surface tension. To address the latter problem we think that the arc length formulation of the 2D water wave problem is the most adapted framework since in this formulation the term with the most derivatives is linear.

Moreover, since the 2D water wave equations with finite depth can be obtained from the 2D water wave equations with infinite depth by replacing ${\rm H}$ by $\k0\,$, we expect that the techniques developed in the present paper and in \cite{D16, DS06} can be generalized and applied to the 2D water wave problem with finite depth and no surface tension as well as 
to the 2D water wave problem with finite depth and with surface tension in the arc length formulation
to derive high-order energy estimates and to establish extended life spans for small solutions analogously to \cite{HIT14}. We think that proving such high-order energy estimates will also be an essential step toward solving the problem of global existence of small localized solutions to the 2D water wave problem with finite depth. 

However, since the 2D water wave problem with finite
depth possesses solitary waves solutions, which are localized traveling waves of permanent form, analogous dispersive estimates as in the case of infinite depth cannot be expected. Hence, in order to
show global existence of small localized solutions to the 2D water wave problem with finite depth one may try to derive appropriate estimates characterizing the solitary wave dynamics of the water wave problem and combine them with the high-order energy estimates. 
\medskip

{\bf Notation}. 
We denote the 
Fourier transform of a function $u \in L^2(\R,\K)$, with $\K=\R$ or $\K=\C$ by
$$\mathcal{F}(u)(k) = \widehat{u}(k) = \frac{1}{2\pi} \int_{\R} u(x) e^{-ikx} dx. $$
Let $H^{s}(\R,\K)$ be
the space of functions mapping from $\R$ into $\K$
for which
the norm
$$ \| u \|_{H^{s}(\R,\K)} = \left(\int_{\R} |\widehat{u}(k)|^2 (1+|k|^2)^{s} 
dk \right)^{1/2} $$ 
is finite. We also write $L^2$ and $H^{s}$ instead of $L^2(\R,\R)$ and $H^{s}(\R,\R)$.
Moreover, we use the space
$ L^p(m)(\R,\K) $ defined by $ u \in L^p(m)(\R,\K) \Leftrightarrow u \sigma^m \in L^p(\R,\K)$, where
$ \sigma(x) = (1+x^2)^{1/2}$.

Furthermore, we write $A
\lesssim B$, if $A \leq C B$ for a constant $C>0$, and $A = \mathcal{O}(B)$, if 
$|A| \lesssim B$. 
\medskip

\textbf{Acknowledgment:} The authors thank the referees for their useful comments.

\section{Long time solutions}
\label{sec2}
In this section, we prove Theorem \ref{Theorem 1 Langzeit loesung}. To address this issue, we will need the following properties of the operator $\k0.$

\begin{lemma}
\label{k0-id-lem}
Let $f,g \in L^{2}$ and $fg \in L^{2}$. Then we have
\begin{align}
\label{k0-id}
\k0 \, (fg)-\k0 \,(f)\,g-f\k0 \,(g)=\k0 \;(\k0 \, (f) \k0 \,(g))\,.
\end{align}
\end{lemma}

\textbf{Proof.}
Considering the symbol of $\k0\,$, we obtain
the assertion of the Lemma due to
\begin{align*}
\tanh(k)-\tanh(m)-\tanh(k-m)=-\tanh(k)\tanh(m)\tanh(k-m)
\end{align*}
for all $m,k \in \R$, which can be directly verified.
\qed

\begin{lemma}
$\kk \partial_x$ is a continuous linear operator from $H^{s+1}$ into $H^{s}$ for any $s \ge 0$ and satisfies 
\begin{align}
\| \kk \partial_x f \|_{H^s} \le \|f\|_{H^{s+1}}
\end{align}
for all $f \in H^{s+1}$.
\end{lemma}

\textbf{Proof.}
The assertion of the Lemma is a consequence of 
\begin{align}
\label{kthk}
|k| \le |\tanh(k)| (1+k^2)^{1/2}
\end{align}
for all $k \in \R$, which can be directly verified.
\qed

\begin{lemma}  \label{commutator1}
Let $j \geq 0$, $q > \frac12$, $r \geq \max\{1+q, j\}$ and $u \in H^{r}$.  Then we have the commutator estimate
\begin{equation} \label{kom1}
\|[K_0^{-1}, u] \partial_x u\|_{H^{j}} \lesssim \|u\|_{H^{1+q}} \|u\|_{H^{j}}\,. 
\end{equation}
\end{lemma}

\textbf{Proof.}
Since $u \partial_x u = \frac{1}{2} \partial_x (u^2)$, the Fourier transform of $[K_0^{-1}, u] \partial_x u$ has the two representations
\begin{align*}
\mathcal{F}({[\kk\,, u] \partial_x u}) (k) &= 
- \int_{\mathbb{R}} \Big( \frac{m}{\tanh(k)} - \frac{m}{\tanh(m)} \Big)  \widehat{u}(k -m) \widehat{u}(m)\, dm \\[2mm]
 &= - \int_{\mathbb{R}} \Big( \frac{k}{2\tanh(k)} - \frac{m}{\tanh(m)} \Big)  \widehat{u}(k -m) \widehat{u}(m)\, dm\,.
\end{align*}
Hence, using Young's inequality for convolutions and the Cauchy-Schwarz inequality, we obtain
\begin{align*}
\| [\kk\,, u] \partial_x u \|_{H^j} 
&\lesssim  \Big( \sup_{k,m \in \mathbb{R}} G(k,m) \Big)\, \|\widehat{u} \|_{L^{2}(1+q)} \|\widehat{u}\|_{L^2(j)}\\
&\lesssim \Big( \sup_{k,m \in \mathbb{R}} G(k,m) \Big) \, \| u \|_{H^{1+q}} \|u\|_{H^j}
\end{align*}
with
\begin{align*}
G(k,m)=\begin{cases}
\Big\vert \dfrac{k}{2\tanh(k)} - \dfrac{m}{\tanh(m)}\Big\vert \, \dfrac{(1+k^2)^{j/2}}{(1+(k-m)^2)^{j/2} (1+m^2)^{1/2}}  & \quad {\rm for }\; |k| \leq 1\,,\\[4mm]
\Big\vert \dfrac{m}{\tanh(k)} - \dfrac{m}{\tanh(m)}\Big\vert \, \dfrac{(1+k^2)^{j/2}}{(1+(k-m)^2)^{j/2} (1+m^2)^{1/2}}  & \quad{\rm for}\; |k| > 1 \,,
\end{cases}
\end{align*}
where $y \mapsto y/\tanh(y)$ is continued by $1$ for $y=0$.

In order to show the boundedness of the supremum we distinguish three cases.
$G(k,m)$ is obviously uniformly bounded for all $(k,m) \in \R^2$ with
$\vert k \vert \le 1$ or $\vert m \vert \le 1$. 
If $\vert k \vert, \vert m \vert \ge 1$ and $\sign(m) = \sign(k)$, we have 
\begin{align*}
\Big\vert \frac{1}{\tanh(k)} - \frac{1}{\tanh(m)} \Big\vert \lesssim e^{-2|k|}+e^{-2|m|}
\end{align*}
such that 
\begin{align*}
G(k,m)
\lesssim
\frac{e^{-2|k|} (1+k^2)^{j/2} \vert m \vert }{(1+(k-m)^2)^{j/2} (1+m^2)^{1/2}} +
\frac{e^{-2|m|} (1+2(k-m)^2+2m^2)^{j/2} \vert m \vert }{(1+(k-m)^2)^{j/2} (1+m^2)^{1/2}} \,.
\end{align*}
Consequently, $G(k,m)$ is also uniformly bounded in this case. Finally,
if $\vert k \vert, \vert m \vert \ge 1$ and $\sign(m) = -\sign(k)$, we have 
\begin{align*}
\Big| \frac{1}{\tanh(k)} - \frac{1}{\tanh(m)}\Big| \lesssim 1
\end{align*}
such that
\begin{align*}
G(k,m)
\lesssim  \frac{ (1+ k^2)^{j/2}  \vert  m \vert}{(1+(|k|+|m|)^2)^{j/2} (1+m^2)^{1/2}} 
\end{align*}
and $G(k,m)$ is uniformly bounded in this case as well. Hence, the supremum is bounded, which implies the assertion of the lemma.
\qed
\medskip

Moreover, we will use the well-known interpolation inequalities
\begin{align}
\label{intpol1}
\Vert \partial_x^{j} f   \partial_x^{\ell} g \Vert_{L^2} 
\lesssim 
\Vert g \Vert_{L^{\infty}} 
\Vert \partial_x^{j+\ell} f \Vert_{L^2} +
\Vert \partial_x^{j+\ell} g \Vert_{L^2} 
\Vert  f \Vert_{L^{\infty}} 
\,,
\end{align}
\begin{align}
\label{intpol2}
\Vert \partial_x^{j} f   \partial_x^{\ell} g \Vert_{L^2} 
\lesssim  
\Vert g \Vert_{H^q} 
\Vert \partial_x^{j+\ell} f \Vert_{L^2} +
\Vert \partial_x^{j+\ell} g \Vert_{L^2} 
\Vert  f \Vert_{H^q} 
\end{align}
for $j,\ell \in \N$, $k \ge 1$, $j+\ell \le k$, $\frac{1}{2}< q \le k$, and $f,g \in H^{k}$ as well as the identity 
\begin{align}
\label{PI}
\int_{\R}f g \partial_x f\, dx = - \frac{1}{2} \int_{\R}f^{2} \partial_x g\,  dx
\end{align}
for $f,g \in H^1$.

\medskip

As motivated in Section \ref{sec1}, we define the energy
\begin{align}
\mathcal{E}_{s}:=\sum_{\ell=0}^{s} E_{\ell}\,,
\end{align}
with
\begin{align}
E_0 :=& \;\frac{1}{2}  \| u\|_{L^2}^2
\end{align}
and
\begin{align}
E_{\ell} :=& \;\frac{1}{2} \| \partial_x^{\ell} u\|_{L^2}^2  
+\sum_{a=1}^{\ell-1} {\binom \ell a} \int_{\mathbb{R}} \kk \partial_x^{\ell}u \kk \partial_x^{a} u \kk \partial_x^{\ell-a+1}u\, dx
\\[2mm]
& \nonumber \;
+ \frac{1}{2} \int_{\mathbb{R}} \kk \partial_x^{\ell}u \kk \partial_x^{\ell}u \kk \partial_xu \,dx
\end{align}
for $\ell >0$.

\begin{remark}
{\rm  One may wonder why we do not write $E_{\ell}$ for $\ell >0$ in the equivalent form
\begin{align*}
E_{\ell} =& \;\frac{1}{2} \| \partial_x^{\ell} u\|_{L^2}^2  
+ \sum_{a=2}^{\ell-1} {\binom {\ell} a} 
\int_{\mathbb{R}}
\kk \partial_x^{\ell} u
\kk \partial_x^a u \kk \partial_x^{\ell-a+1} u \,dx 
\\[2mm] &
+ \frac{2\ell+1}{2}
\int_{\mathbb{R}}
\kk \partial_x^{\ell} u
\kk \partial_x^{\ell} u \kk \partial_x u\, dx\,,
\end{align*}
but it turns out that the form \eqref{Eell-def} is more convenient for our
calculations below.}
\end{remark}
The Cauchy-Schwarz inequality and \eqref{intpol1} directly imply
\begin{lemma}
\label{lem Energie aequi}
For $\ell\ge 1$, we have
\begin{align}
\label{lin  Energie aequi}
E_{\ell}=\frac{1}{2} \Vert \partial_x^{\ell} u \Vert_{L^2}^2 
+ \mathcal{O}( \Vert \kk \partial_x u \Vert_{L^\infty} )\Vert \kk \partial_x^{\ell} u \Vert_{L^2}^2\,.
\end{align}
\end{lemma}

Now, we would like to estimate the time derivative of $\mathcal{E}_{s}$ for any sufficiently regular solution of \eqref{DGL}. We obtain

\begin{lemma}
\label{E0-lem}
For $\ell=0$, we have
\begin{align}
\label{lin dt E, l=0}
\frac{d}{dt} E_0=0\,.
\end{align}
\end{lemma}
\textbf{Proof.}
Due to  \eqref{DGL}, the skew symmetry of $\k0\,$ and \eqref{PI}, we get
\begin{align*}
\frac{d}{dt} E_{0} =& 
\int_{\mathbb{R}} u \k0 u\, dx
-\int_{\mathbb{R}}  u^{2} \,\partial_x u\, dx \,
=\,0\,.
\end{align*}
\qed

\begin{lemma}
For $\ell\ge 1$, we have
\begin{align}
\label{lin dt E, l>0}
\frac{d}{dt} E_{\ell} \lesssim \Vert u \Vert_{H^2}^2 \Vert u\Vert _{H^{\ell}}^2\,.
\end{align}
\end{lemma}

\textbf{Proof.}
Using \eqref{DGL}, we get 
\begin{align*}
\frac{d}{dt} E_{\ell} =& 
\int_{\mathbb{R}} \partial_x^{\ell} u \, \partial_x^{\ell}\!\k0 u\, dx
-\int_{\mathbb{R}} \,\partial_x^{\ell} u \,\partial_x^{\ell} (u \,\partial_x u)\, dx \\
&
+\sum_{a=1}^{\ell-1} {\binom \ell a} \int_{\mathbb{R}}   \partial_x^{\ell}u  \kk \partial_x^{a} u  \kk \partial_x^{\ell-a+1}u \,dx
+ \frac{1}{2} \int_{\mathbb{R}}  \partial_x^{\ell}u \kk \partial_x^{\ell}u  \kk \partial_xu \,dx
\\
&
+\sum_{a=1}^{\ell-1} {\binom \ell a} \int_{\mathbb{R}}  \kk\partial_x^{\ell}u \, \partial_x^{a} u \kk \partial_x^{\ell-a+1}u\, dx
+ \frac{1}{2} \int_{\mathbb{R}}  \kk\partial_x^{\ell}u \, \partial_x^{\ell}u \kk \partial_xu \,dx
\\
&
+\sum_{a=1}^{\ell-1} {\binom \ell a} \int_{\mathbb{R}}  \kk\partial_x^{\ell}u  \kk\partial_x^{a} u \,  \partial_x^{\ell-a+1}u \,dx
+ \frac{1}{2} \int_{\mathbb{R}}  \kk\partial_x^{\ell}u \kk \partial_x^{\ell}u \, \partial_xu \,dx
\\
&
-\sum_{a=1}^{\ell-1} {\binom \ell a} \int_{\mathbb{R}} \kk \partial_x^{\ell}(u \, \partial_x u) \kk \partial_x^{a} u \kk \partial_x^{\ell-a+1}u \,dx
\\
&
-\sum_{a=1}^{\ell-1} {\binom \ell a} \int_{\mathbb{R}} \kk \partial_x^{\ell}u \kk \partial_x^{a} (u\, \partial_x u)\kk \partial_x^{\ell-a+1}u \,dx
\\
&
-\sum_{a=1}^{\ell-1} {\binom \ell a} \int_{\mathbb{R}} \kk \partial_x^{\ell}u \kk \partial_x^{a} u \kk \partial_x^{\ell-a+1}(u\, \partial_x u) \,dx
\\&
-\int_{\mathbb{R}} \kk \partial_x^{\ell}(u\, \partial_x u) \kk \partial_x^{\ell}u \kk \partial_xu \,dx
- \frac{1}{2} \int_{\mathbb{R}} \kk \partial_x^{\ell}u \kk \partial_x^{\ell}u \kk \partial_x (u \,\partial_x u) \,dx\,.
\end{align*}
Due to the skew symmetry of $\k0\,$, the first integral equals zero.
Because of \eqref{partintf}, \eqref{k0-id} and the skew symmetry of $\k0 \,$ all integrals with cubic integrands cancel. We recall that this cancellation is a consequence of the identities
\eqref{Nrel} and \eqref{Ndef1} of the normal-form transform $N$ and that $E_{\ell}$ was constructed by including this normal-form transform in order to obtain that the right-hand side of the evolution equation of $E_{\ell}$ consists only of quartic terms.

Hence, we have
\begin{align*}
\frac{d}{dt} E_{\ell} =& 
-\sum_{a=1}^{\ell} {\binom \ell a} \int_{\mathbb{R}} \kk \partial_x^{\ell}(u \, \partial_x u) \kk \partial_x^{a} u \kk \partial_x^{\ell-a+1}u \,dx
\\
&
-\sum_{a=1}^{\ell-1} {\binom \ell a} \int_{\mathbb{R}} \kk \partial_x^{\ell}u \kk \partial_x^{a} (u\, \partial_x u)\kk \partial_x^{\ell-a+1}u \,dx
\\
&
-\sum_{a=1}^{\ell-1} {\binom \ell a} \int_{\mathbb{R}} \kk \partial_x^{\ell}u \kk \partial_x^{a} u \kk \partial_x^{\ell-a+1}(u\, \partial_x u) \,dx
\\&
- \frac{1}{2} \int_{\mathbb{R}} \kk \partial_x^{\ell}u \kk \partial_x^{\ell}u \kk \partial_x (u \,\partial_x u) \,dx\,.
\end{align*}
Since ${\binom \ell a} = {\binom \ell {\ell-a} }$, we obtain by integration by parts
\begin{align*}
\frac{d}{dt} E_{\ell}= \sum_{a=1}^{\ell} {\binom \ell a} I_{a}\,,
\end{align*}
with
\begin{align*}
I_{a}=& - \int_{\mathbb{R}}\kk \partial_x^{\ell}(u\, \partial_x u)\kk \partial_x^{a}u \kk \partial_x^{\ell- a +1} u\,  dx 
\\
&
+ \int_{\mathbb{R}} \kk \partial_x^{\ell+1} u \kk \partial_x^{\ell-a}(u\, \partial_x u) \kk \partial_x^{a} u \,dx
\\
=& - \int_{\mathbb{R}} \partial_x^{\ell}(u \kk \partial_x u )\kk \partial_x^{a}u \kk \partial_x^{\ell- a +1} u\,  dx 
\\
&
+ \int_{\mathbb{R}} \kk \partial_x^{\ell+1} u \, \partial_x^{\ell-a}(u \kk \partial_x u) \kk \partial_x^{a} u \,dx
\\
&
- \int_{\mathbb{R}} \partial_x^{\ell}[\kk, u] \partial_x u
\kk \partial_x^{a}u \kk \partial_x^{\ell- a +1} u\, dx
\\
&
+ \int_{\mathbb{R}} \kk \partial_x^{\ell+1} u \, \partial_x^{\ell-a}[\kk, u] \partial_x u  \kk \partial_x^{a} u\, dx
\\
=:& \sum_{k=1}^{4} I_{a,k}\,.
\end{align*}
With the help of Leibniz's rule we get
\begin{align*}
I_{a,1}+ I_{a,2}=& - 
\int_{\mathbb{R}} u \kk \partial_x^{\ell+1} u
\kk \partial_x^{a}u \kk \partial_x^{\ell- a +1} u 
\, dx
\\
&
- \sum_{i=1}^{\ell}
{ \binom  \ell i}
\int_{\mathbb{R}} \partial_x^{i} u \kk \partial_x^{\ell-i+1} u
\kk \partial_x^{a}u \kk \partial_x^{\ell- a +1} u 
\, dx
\\
&
+\int_{\mathbb{R}} \kk \partial_x^{\ell+1} u\,
u \kk  \partial_x^{\ell-a +1}u \kk \partial_x^{a} u \,dx
\\
&
+\sum_{j=1}^{\ell-a} { \binom  {\ell-a} j}
\int_{\mathbb{R}} \kk \partial_x^{\ell+1} u\,
\partial_x^{j}u \kk  \partial_x^{\ell-a -j+1}u \kk \partial_x^{a} u \,dx
\\
=& 
- \sum_{i=1}^{\ell}
{ \binom  \ell i}
\int_{\mathbb{R}} \partial_x^{i} u \kk \partial_x^{\ell-i+1} u
\kk \partial_x^{a}u \kk \partial_x^{\ell- a +1} u 
\, dx
\\
&
+\sum_{j=1}^{\ell-a} { \binom  {\ell-a} j}
\int_{\mathbb{R}} \kk \partial_x^{\ell+1} u\,
\partial_x^{j}u \kk  \partial_x^{\ell-a -j+1}u \kk \partial_x^{a} u \,dx
\,.
\end{align*}
Using the interpolation inequality \eqref{intpol2} and the Cauchy-Schwarz inequality yields
\begin{align*}
\sum_{i=1}^{\ell}
{ \binom  \ell i}
\int_{\mathbb{R}} \partial_x^{i} u \kk \partial_x^{\ell-i+1} u
\kk \partial_x^{a}u \kk \partial_x^{\ell- a +1} u 
\, dx 
\lesssim &\, \Vert u \Vert_{H^2}^2 \Vert u\Vert _{H^{\ell}}^2 
\,.
\end{align*}
The remaining integrals
\begin{align*}
\int_{\mathbb{R}} \kk \partial_x^{\ell+1} u\,
\partial_x^{j}u \kk  \partial_x^{\ell-a -j+1}u \kk \partial_x^{a} u \,dx
\end{align*}
can be rewritten by applying a finite, $\ell$-dependent number of integrations by parts into a sum of integrals of the form
\begin{align*}
\int_{\mathbb{R}} \partial_x^{\alpha} u\,
\kk \partial_x^{\beta}u \kk  \partial_x^{\gamma}u \kk \partial_x^{\delta} u \,dx
\end{align*}
with
\begin{align*}
1 \leq \alpha,\beta,\gamma,\delta  \quad \mathrm{and}
\quad \alpha+\beta = \gamma+\delta = \ell +1
\end{align*}
such that we can apply again \eqref{intpol2} and the Cauchy-Schwarz inequality to obtain
\begin{align*}
\sum_{j=1}^{\ell-a} { \binom  {\ell-a} j}
\int_{\mathbb{R}} \kk \partial_x^{\ell+1} u\,
\partial_x^{j}u \kk  \partial_x^{\ell-a -j+1}u \kk \partial_x^{a} u \,dx
\lesssim &\, \Vert u \Vert_{H^2}^2 \Vert u\Vert _{H^{\ell}}^2 
\,.
\end{align*}
Hence, we have shown
\begin{align*}
I_{a,1}+ I_{a,2}
\lesssim &\, \Vert u \Vert_{H^2}^2 \Vert u\Vert _{H^{\ell}}^2 
\,.
\end{align*}
Finally, using the Cauchy-Schwarz inequality, \eqref{kom1} and \eqref{intpol2}, we get
$$
I_{a,3} \lesssim \Vert u \Vert_{H^2}^2 \Vert u\Vert _{H^{\ell}}^2\,,
$$
and with the aid of integration by parts, the Cauchy-Schwarz inequality, \eqref{kom1},  \eqref{intpol2} and \eqref{PI}, we obtain
$$
I_{a,4} \lesssim \Vert u \Vert_{H^2}^2 \Vert u\Vert _{H^{\ell}}^2\,.
$$
\qed
\medskip

Now, combining the estimates \eqref{lin dt E, l=0}, \eqref{lin dt E, l>0} and \eqref{lin  Energie aequi}, we get  
\begin{equation}
\frac{d}{dt} \mathcal{E}_{s} \lesssim  \varepsilon^2 \mathcal{E}_{s}
\end{equation}
for any solution $u \in C(I, H^s)\, \cap\, C^1(I, H^{s-1}) $, where $I \subset \R$ and $s \ge 2$, of \eqref{DGL} with $\|u\|_{H^{2}} \leq \eps$.  
Because of the local existence results for quasi-linear symmetric hyperbolic systems 
from \cite{Kat75} and Gronwall's inequality, we obtain  
the $\mathcal{O}(1)$-boundedness of $\mathcal{E}_{s}$ and therefore of $\|u\|_{H^{s}}$ for all $t \in I_{\eps}$, which proves Theorem \ref{Theorem 1 Langzeit loesung}.
\qed

\section{The derivation of the NLS approximation}
\label{sec3}
In this section, we derive the NLS equation as an approximation equation for the quasilinear dispersive equation \eqref{DGL}.
In doing so, we make the ansatz 
\begin{equation}
\label{ansatz}
u= \varepsilon \widetilde{\psi} =
\varepsilon \widetilde{\psi}_1+\varepsilon \widetilde{\psi}_{-1}
+\varepsilon^2 \widetilde{\psi}_0
+\varepsilon^2 \widetilde{\psi}_2+\varepsilon^2 \widetilde{\psi}_{-2}\,,
\end{equation}
with
\begin{equation*}
\widetilde{\psi}_j(x,t) = \widetilde{A}_j(\varepsilon(x-c_g t), \varepsilon^2 t)\, \textbf{E}^{j}
\end{equation*}
and $\widetilde{A}_{-j} = \overline{\widetilde{A}_j}$
for $j \in \{0,1,2\}$, where
$0< \varepsilon \ll 1$, $ k_0 >0$, $\omega_0= \tanh(k_0)$, $c_g= \tanh'(k_0) = {\sech}^2(k_0)$, and $\textbf{E}= e^{i(k_0x-\omega_0t)}$.
\begin{remark}{ Our ansatz leads to waves moving to the right.
For waves moving to the left one has to replace in the above ansatz
$\omega_0$ by $-\omega_0$
and $c_g$ by $-c_g$.
}
\end{remark}

We insert our ansatz \eqref{ansatz} in equation \eqref{DGL}. Then we expand all terms of the form $\k0 \widetilde{\psi}_j$ by using the Taylor series of the hyperbolic tangent around $k=jk_0$. (For more details compare Lemma 25 in \cite{SW10}, for example.) After that we equate the coefficients in front of the $\varepsilon^m \textbf{E}^{j}$ to zero.
In detail, we get for 
\begin{align*}
& (m,j)=(1,1) : && i\omega_0\widetilde{A}_1= i \tanh(k_0)\widetilde{A}_1\,, \\
& (m,j)=(2,1) : && c_g \partial_X \widetilde{A}_1= {\sech}^2(k_0) \partial_X \widetilde{A}_1\,, \\
& (m,j)=(2,2) : && i(-2\omega_0+\tanh(2k_0))\widetilde{A}_2= ik_0 (\widetilde{A}_1)^2\,,\\
& (m,j)=(3,0) : && (-c_g + {\sech}^2(0)) \partial_X \widetilde{A}_0  = \partial_X(\widetilde{A}_1 \widetilde{A}_{-1})\,, \\
& (m,j)=(3,1) : && \partial_T \widetilde{A}_1= -i\, \tanh(k_0)\, {\sech}^2(k_0)\, \partial_X^2 \widetilde{A}_1 +ik_0(\widetilde{A}_0\widetilde{A}_1+\widetilde{A}_{-1}\widetilde{A}_2)\,,
\end{align*}
where $X=\varepsilon(x-c_g t)$ and $T= \eps^2t$.

The equations for $(m,j)=(1,1)$ and $(m,j)=(2,1)$ are satisfied due to the definitions of $\omega_0$ and $c_g$.
Since for $k_0 \neq 0$ and all integers $j \geq 2$ the non-resonance conditions
\begin{equation}
\label{NR1a}
\tanh(jk_0)\neq j \tanh(k_0)\,,
\end{equation}
\begin{equation}
\label{NR1b}
\tanh'(k_0) \neq \tanh'(0)
\end{equation}
hold, we can choose $\widetilde{A}_0$ and $\widetilde{A}_2$ depending on $\widetilde{A}_1$, such that the equations for $(m,j)=(2,2)$ and $(m,j)=(3,0)$ are satisfied and  
the equation for $(j,m)=(3,1) $ becomes the NLS equation
\begin{equation}
\label{NLS-Eq}
\partial_T \widetilde{A}_1= i \nu_1 \partial_X^2 \widetilde{A}_1 +i \nu_2 |\widetilde{A}_1|^2 \widetilde{A}_1\,, 
\end{equation}
with 
\begin{align*}
&\nu_1= \frac{1}{2} \tanh''(k_0)= - \tanh(k_0) {\sech}^2(k_0)\,, 
\\
& \nu_2= k_0 
\left( 
\frac{k_0}{\tanh(2k_0)-2\tanh(k_0)} +\frac{1}{{\tanh}^2(k_0)}
\right).
\end{align*}

To prove the approximation property of the NLS equation \eqref{NLS-Eq} it will be helpful to make the residual 
\begin{equation}
\res (\veps \widetilde{\psi})= -\partial_t (\veps \widetilde{\psi})+ \k0 \,(\veps \widetilde{\psi})- \veps \widetilde{\psi} \partial_x (\veps \widetilde{\psi})\,,
\end{equation}
which contains all terms that do not cancel after inserting ansatz 
\eqref{ansatz} into system \eqref{DGL}, smaller in any Sobolev norm $\| \cdot \|_{H^s}$ with $s \geq 0$ by proceeding analogously as in Section 2 of \cite{DSW12} and replacing $\veps \widetilde{\psi}$ by a new approximation $\eps\psi$ of the form
\begin{equation}
\label{ans-high}
\varepsilon \psi =\sum_{|j| \le 5} \sum_{\beta(j,n) \le 5} \varepsilon^{\beta(j,n)} \psi_{j}^{n}\,,
\end{equation}
where $j \in \mathbb{Z}$, $n \in \N_0$, 
\begin{align}
\label{ext-ans1}
&\beta(j,n)=1+ \vert \vert j \vert -1 \vert +n\,, \\
&\psi_{j}^{n}(x,t)= A_{j}^{n}(\varepsilon(x- c_g t), \varepsilon^2 t)\textbf{E}^{j}\,,
\label{ext-ans2}
\end{align}
${A}_{-j}^{n} = \overline{{A}_j^{n}}$, and the functions
${\psi}_{j}^{n}$ have the compact support 
\begin{align}
&\{k \in \R : |k- jk_0| \le \delta < k_0/20\}
\end{align}
in Fourier space, for sufficiently small $\eps >0$.
For later purposes we fix $\delta \in (0, k_0/20)$ such that 
\begin{equation}
\label{delta}
\vert \tanh(k) -2 \tanh(jk_0)-\tanh(k-2jk_0) \vert \ge C>0
\end{equation}
for a constant $C=C(k_0)$, which is possible due to \eqref{NR1a}.

This new approximation is constructed in the following way. First, the previous approximation $\veps \widetilde{\psi}$ is extended by higher order correction terms such that the resulting approximation, which we denote by $\veps \widetilde{\psi}_{ext}$, has the form \eqref{ans-high}-\eqref{ext-ans2} with $\psi$, $\psi_{j}^{n}$ and ${A}_{j}^{n}$ replaced by $\widetilde{\psi}_{ext}$, $\widetilde{\psi}_{j}^{n}$ and $\widetilde{A}_{j}^{n}$, where 
$\widetilde{A}_{j}^{0}= \widetilde{A}_{j}$ and the higher order correctors
$\widetilde{A}_{j}^{n}$, $n > 0$, can be computed by a similar procedure as the functions $\widetilde{A}_{j}$.  More precisely, inserting $\veps \widetilde{\psi}_{ext}$ into \eqref{DGL} and equating the coefficients in front of the $\varepsilon^{\beta(j,n)} \textbf{E}^{j}$ to zero yields 
a system of algebraic equations and inhomogeneous linear Schr\"odinger equations that can be solved recursively. Due to the non-resonance conditions \eqref{NR1a}-\eqref{NR1b} the functions $\widetilde{A}_{j}^{n}$ with $j \neq \pm 1$ are uniquely determined by the algebraic equations. The functions $\widetilde{A}_{\pm 1}^{n}$ satisfy the inhomogeneous linear Schr\"odinger equations. Moreover, since the functions $\widetilde{A}_{\pm 1}^4$ do not appear in the equations for any other $\widetilde{A}_{j}^{n}$, we can set  $\widetilde{A}_{\pm 1}^4=0$.

Secondly, by multiplying the Fourier transform of each function $\widetilde{\psi}_{j}^{n}$ by a suitable cut-off function, we obtain our final approximation $ \varepsilon \psi$. Since the Fourier transform of the functions $\widetilde{\psi}_{j}^{n}$ is strongly concentrated around the wave number $jk_0$ if $\widetilde{A}_{j}^{n}$ is sufficiently regular, the approximation is only changed slightly by the second modification, but this action will give us a simpler control
of the error and makes the approximation an analytic function.

Furthermore, we define
\begin{align}
&\psi_{\pm1}:=\psi_{\pm1}^{0}\,, \\
&\psi_c:=\psi_{-1}+\psi_{1}\,, \\
&\psi_{s}:= \eps^{-1}(\psi-\psi_c)
\end{align}
and get the following estimates for the modified residual. 
\begin{lemma}
\label{Res-abschaetzung}
Let  $s_A \ge 7$,  $\widetilde{A}_1 \in C([0,T_0],H^{s_A}(\mathbb{R},\mathbb{C}))$ be a solution of the NLS equation \eqref{NLS-Eq} with
\begin{align*}
\sup_{T \in [0, T_0]} \Vert \widetilde{A}_1(T) \Vert_{H^{s_A}} \le C_A\,,
\end{align*}
and $\delta$ be chosen as above.
Then for all $s \geq 0$ there exist $ C_{\rm Res}, C_{\psi}, \eps_0>0  $ depending on  $C_A$, $k_0$ and $\delta$, where $\eps_0 < \delta$,  such that for all
$ \eps \in (0,\eps_0) $ the approximation $\eps \psi$ satisfies
\begin{eqnarray} \label{RES1}
\sup_{t \in [0,T_0/\eps^2]} \|  {\rm Res}(\eps \psi)
\|_{H^s}
& \leq & C_{\rm Res}\, \eps^{{11/2}}, \\ \label{RES2}
\sup_{t \in [0,T_0/\eps^2]} \|\eps \psi - (\eps {\widetilde{\psi}_{1}} + \eps {\widetilde{\psi}_{-1}})
\|_{H^{{s_A}}}
& \leq & {C_{\psi}}\, \eps^{3/2},
\\
\label{RES3}
\sup_{t \in [0,T_0/\eps^2]} (\|\widehat{\psi}_{\pm1}\|_{L^1({s})}+\|\widehat{\psi}_s\|_{L^1({s})})
&{\leq} & {C_{\psi}}\,. 
\end{eqnarray}
\end{lemma} 

\textbf{Proof.}
The first extended approximation $\veps \widetilde{\psi}_{ext}$ is constructed in a way that formally we have ${\rm Res}(\veps \widetilde{\psi}) = \OO(\epsilon^{6})$ and $\veps\widetilde{\psi} - (\veps\widetilde{\psi}_{1}+\veps\widetilde{\psi}_{-1}) = \OO(\epsilon^{2})$ on the time interval $[0,T_0/\epsilon^2]$ if $\widetilde{A}_1$
is a solution of the NLS equation \eqref{NLS-Eq} for $T \in [0,T_0]$. 

It can be shown exactly as in the proof of Theorem 2.5 in \cite{DSW12} that  
 $ \widetilde{A}_1 \in C([0,T_0],H^{s_A}) $ with $s_A \geq 5$ implies
$ \widetilde{A}_{j}^n \in C([0,T_0],H^{s_A-{{n}}}) $ if $j \neq \pm1 $ and
$ \widetilde{A}_{\pm 1}^n \in C([0,T_0],H^{s_A-n-2})$ for $n \in \{1,2,3\}$, where the
respective Sobolev norms are uniformly bounded by the $H^{s_A}$-norm of  $ \widetilde{A}_1$.
Therefore, by taking into account that
$\| f(\epsilon\, \cdot)\|_{L^2} = \epsilon^{-1/2} \| f \|_{L^2}$, we obtain
 estimates of the form \eqref{RES1} and \eqref{RES2} with $\psi$ replaced by $ \widetilde{\psi}_{ext}$ and $H^{s}, H^{s_A}$ replaced by $L^2$ if we have $ \widetilde{A}_1 \in C([0,T_0],H^{s_A}) $ with $s_A \geq 7$ (since two additional 
spatial derivatives of $\widetilde{A}_1$ are needed to bound ${\rm Res} (\veps \widetilde{\psi}_{ext})$ in $L^2$).

Since the Fourier transform of the final approximation $\eps \psi$ has a compact support whose size depends on $k_0$, there exists a $C=C(k_0) >0$ such that $\|\psi\|_{H^s} \leq C \|\psi\|_{L^2}$ and 
$\|\widehat{\psi}\|_{L^1({s})} \leq C \|\widehat{\psi} \|_{L^1}$ for all $s \geq 0$.
Hence, by using the above $L^2$-estimates for $\veps \widetilde{\psi}_{ext}$ as well as the estimate
\begin{equation}  \label{cut} 
 \| (\chi_{[-\delta,\delta]}-1)\, \eps^{-1} 
\widehat{f} ( \eps^{-1} \cdot) \|_{L^2(m)}  \leq C(\delta)\,
\eps^{m+{M}-1/2} \| f \|_{H^{m+{M}}}
\end{equation}
for all $M,m \geq 0$, where $\chi_{[-\delta,\delta]}$ is the characteristic function on $[-\delta,\delta]$, for $f=\widetilde{A}_{j}^n$ for each $\widetilde{A}_{j}^n$ with $m=0$, $M=M(j,n)$ determined by the maximal Sobolev regularity of the respective $\widetilde{A}_{j}^n$ and $\delta$ as above, we obtain \eqref{RES1} and
\begin{equation}  \label{RES4}
\sup_{t \in [0,T_0/\eps^2]} \|\eps (\psi-\psi_c)
\|_{H^{{s_A}}}
 \leq  {C_{\psi}}\, \eps^{3/2}
\end{equation}
if we have $ s_A  \geq 7$, which yields $\beta(j,n)+ M(j,n) \geq 6$.  
By combining \eqref{RES4} and \eqref{cut} for $f=\eps {\widetilde{\psi}_{1}} + \eps {\widetilde{\psi}_{-1}}$, $m=s_A$, $M=0$ and $\delta$ as above, we obtain \eqref{RES2}.

Finally, since $\| \epsilon^{-1} \widehat{f}(\epsilon^{-1}  \cdot)\|_{L^1} = \| \widehat{f} \|_{L^1}$,
estimate \eqref{RES3} follows by construction of ${\psi}_{\pm 1}$ and ${\psi}_s$.
\qed

\begin{remark} \label{remneucam2}
{\rm 
The bound \eqref{RES3} will be  used for instance to estimate 
$$ \| \psi_{j}^{n}f \|_{H^{s}} \leq C \| \psi_{j}^{n}  \|_{C^{s}_b}\| f 
\|_{H^{s}} \leq  C \| \widehat{{\psi}}_{j}^{n} \|_{L^1(s)}\| f
\|_{H^{s}} $$
without loss of powers in $ \eps $ as it would be the case with $\| \widehat{{\psi}}_{j}^{n}  \|_{L^2(s)}$.}
\end{remark}

Moreover, by an analogous argumentation as in the proof of Lemma 3.3 in \cite{DSW12} we obtain the fact that $\partial_t \psi_{\pm 1}$ can be approximated by ${\k0 \psi}_{\pm 1}$. More precisely, we get

\begin{lemma}
For all $s>0$ there exists a constant $C>0$ depending on $ \Vert \widetilde{A}_1 \Vert_{H^{3}}$ and $k_0$ such that
\begin{equation}
\label{d/dt Psi}
\Vert \partial_t \widehat{\psi}_{\pm 1} - {\k0 \widehat{\psi}}_{\pm 1}  \Vert_{L^1(s)}\le C \varepsilon^2.
\end{equation}
\end{lemma}

\section{The error estimates}
\label{sec4}
Now, we write a solution $u$ of \eqref{DGL} as the sum of approximation and error.
To avoid problems arising from the resonances at $k=\pm k_0$, we rescale the error
with the help of the weight function
\begin{equation}
\widehat{\vartheta}(k) = \left\{ \begin{array}{cc} 1  &\quad {\rm for }\; |k| > \delta\, ,\\[1mm]
\epsilon + (1-\epsilon) |k|/\delta  &\quad {\rm for }\; |k| \leq \delta \, ,
\end{array} \right.
\end{equation}
where $\delta$ is chosen as above and $ \eps \in (0,\eps_0) $, with $\eps_0$ as in Lemma \ref{Res-abschaetzung}. That means, we write
\begin{align}
u= \eps \psi +\eps^{5/2} \vartheta R\,,
\end{align}
where $\vartheta R$ is defined by $\widehat{\vartheta R} = \widehat{\vartheta} \widehat{R}$. By this choice $\widehat{\vartheta R}(k) 
$ is small at the wave numbers close to zero reflecting 
the fact that the nonlinearity of \eqref{DGL} vanishes at $k = 0$. 

Inserting this ansatz into \eqref{DGL} leads to
\begin{equation}
\label{d/dt R}
\partial_t R= \k0R - \eps \vartheta^{-1} \partial_x(\psi \vartheta R)
-\frac{1}{2} \varepsilon^{5/2}\vartheta^{-1} \partial_x (\vartheta R)^2 +\varepsilon^{-5/2} \vartheta^{-1}\res(\eps \psi)\,,
\end{equation}  
where the operator $\vartheta^{-1}$ is defined by its symbol $\widehat{\vartheta^{-1}}(k)= \widehat{\vartheta}^{-1}(k)               =(\widehat{\vartheta}(k))^{-1}$.

Due to the structure of the nonlinear terms in the error equation \eqref{d/dt R}, the
size of the Fourier transform of these terms depends on whether $k$ is close 
to zero or not. In order to separate the behavior in these two
regions more clearly, we define projection operators $P_{0,\alpha}$ and $P_{\alpha,\infty}$ for $\alpha >0$ by the Fourier
multipliers 
\begin{align}
&\widehat{P}_{0,\alpha}(k) = \chi_{[-\alpha, \alpha]}(k)\,,\\[2mm]
&\widehat{P}_{\alpha,\infty}(k)= (1- \chi_{[-\alpha, \alpha]})(k)\,,
\end{align}
where $\chi_{[-\alpha, \alpha]}$ is the characteristic
function on $[-\alpha, \alpha]$.

As motivated in Section \ref{sec1}, we define the energy
\begin{equation}
\tilde{\mathcal{E}}_{s} =\sum_{\ell=0}^{s} \tilde{E}_{\ell}\,,
\end{equation}
\medskip
\begin{equation}
\tilde{E}_{\ell} = \left\{ \begin{array}{ll}  \Vert \check{R}\Vert_{L^2}^2 &\quad {\rm for }\;\ell=0\,,
\\[3mm] 
\dfrac{1}{2} \| \partial_x^{\ell} R\|_{L^2}^2 
+ \veps \displaystyle \int_{\mathbb{R}}  \partial_x^{\ell}R\,  \partial_x^{\ell} N(\psi_c, R)\, dx
&\quad {\rm for }\;\ell >0 \end{array} \right.
\end{equation}
\medskip
with
\begin{equation}
\label{ch R}
\check{R}=R+ \eps N(\psi_c, R)+ \eps^2 \mathcal{T}(\psi_c,\psi_c,R)\,, 
\end{equation}
\begin{equation}
N(\psi_c, R) = - \vartheta^{-1} \kk \partial_x( \kk \psi_c \kk \vartheta P_{\varepsilon, \infty} R)\,,
\end{equation}
\begin{equation}
\widehat{\mathcal{T}}(\psi_c,\psi_c,R)(k)= \sum_{j=\pm1}
\widehat{\mathcal{T}}_j(\psi_j,\psi_j,R)(k)
\,,
\end{equation}
\begin{equation}
\widehat{\mathcal{T}}_j(\psi_j,\psi_j,R)(k)= 
\int_{\mathbb{R}}\int_{\mathbb{R}}
\widehat{t}_j(k)\, \widehat{\psi}_j(k-m)\, \widehat{\psi}_j(m-n)\,\hat{R}(n)\, dn\, dm
\,,
\end{equation}
\begin{align}
\widehat{t}_j(k)= &-\frac{k(k-jk_0)\,\widehat{\vartheta}(k-2jk_0)\,\chi_{[-\delta, \delta]}(k)}{\widehat{\vartheta}(k) 
\tanh(k) \tanh(jk_0) \tanh(k-jk_0)} \\[2mm] 
\nonumber & \times (\tanh(k) -2 \tanh(jk_0)-\tanh(k-2jk_0))^{-1} \,,
\end{align}
\medskip
where $s=s_A\geq 7$, in order to control the error.

To perform our energy estimates we will need the following lemmas.

\begin{lemma}
\label{Nlem}
The operator $N$ has the following properties:\\[2mm]
{\bf a)}
$f \mapsto N(\psi_c,f)$ defines a continuous linear map from $H^1(\R,\R)$ into $L^2(\R,\R)$, and there exists a constant $C=C(\psi_c)>0$, such that for all $f\in H^1(\R,\R)$ and all $g\in H^2(\R,\R)$ we have
\begin{align}
\label{N0}
\Vert N(\psi_c, f) \Vert_{L^{2}} \leq C \varepsilon^{-1} \Vert f \Vert_{H^1}\,,
\end{align}
\begin{align}
\label{PN}
\Vert P_{\delta, \infty}N(\psi_c, f) \Vert_{L^{2}} \leq C \Vert f \Vert_{H^1}\,,
\end{align}\begin{align}
\label{d/dxN}
\Vert \partial_x N(\psi_c, g) \Vert_{L^{2}} \leq C \Vert g \Vert_{H^2}\,.
\end{align}
{\bf b)} 
For all $f\in H^1(\R,\R)$ we have
\begin{align}
\label{N aufloesen}
\vartheta N(\psi_c,f) =  \partial_x ( \kk \psi_c  f) + Q(\psi_c, f)
\end{align}
with
\begin{align}
\label{Q0}
\Vert Q(\psi_c, f) \Vert_{H^{s}}= \mathcal{O}(\Vert f \Vert_{L^2})
\end{align}
for all $s \ge 0$.
\\
{\bf c)}
For all $f \in H^1(\R,\R)$ we have
\begin{align}
\label{N-wahl}
-\k0N(\psi_c,R)
+N(\k0 \psi_c,R)
+N(\psi_c,\k0R) =\vartheta^{-1} \partial_x(\psi_c \vartheta P_{\varepsilon, \infty} R)\,.
\end{align}
{\bf d)}
For all $f \in L^2(\R,\R)$ we have
\begin{align}
\label{P_0 N}
P_{0, \delta}N(\psi_c, P_{0, \delta}f)=0\,.
\end{align}
{\bf e)}\;
For all $f,g \in H^1(\R,\R)$ we have
\begin{equation} \label{partN}
\int_{\R} f\, \vartheta N(\psi_c,g)\,dx = - \int_{\R} g\, \vartheta N(\psi_c,f) \,dx + \int_{\R} S(\partial_x \psi_c,f)\, g\,dx + \int_{\R} Z(\psi_c,f,g)\,dx\,,
\end{equation}
where
\begin{equation*}
S(\partial_x \psi_c,f)= \kk \partial_x \psi_cf \,,
\end{equation*}
\begin{equation*}
Z(\psi_c,f,g) = f\,Q(\psi_c,g) + g\,Q(\psi_c,f)\,.
\end{equation*}
\end{lemma}

\textbf{Proof.}
In Fourier space, we have 
\begin{align}
\label{NFourier}
\widehat{N}(\psi_c,f)(k)= \int_{\R}\widehat{n}(k, k-m,m)\, \widehat{\psi}_c(k-m)\,\widehat{f}(m)\, dm
\end{align}
with
\begin{align*}
\widehat{n}(k,k-m,m)=- {\widehat{\vartheta}}^{-1}(k)\, \widehat{K}_0^{-1}(k)\, ik\,  \widehat{K}_0^{-1}(k-m)\,  \chi_c(k-m)\, \widehat{K}_0^{-1}(m)\, \widehat{\vartheta}(m)\, \widehat{P}_{\varepsilon, \infty}(m)\,,
\end{align*}
where $\chi_c= \chi_{\,\supp(\psi_c)}$.
Now, we estimate the kernel $\widehat{n}$. We have
\begin{align*}
|\widehat{K}_0^{-1}(m)\, \widehat{\vartheta}(m)\, \widehat{P}_{\varepsilon, \infty}(m)|
= 
\left\{ \begin{array}{ll}
 0  & \quad {\rm for }\; 0 < |m| \le \varepsilon\,, \\[2mm]
\dfrac{\varepsilon}{|\tanh(m)|} + \dfrac{(1 - \varepsilon )|m|}{\delta |\tanh(m)|}  & \quad {\rm for }\; \varepsilon \le |m| \le \delta\, , \\[4mm]
  \dfrac{1}{|\tanh(m)|} & \quad {\rm for }\;  |m| \ge \delta \,.
  \end{array} \right.
\end{align*}
Exploiting the monotonicity properties of $m \mapsto 1/|\tanh(m)|$ and $ m \mapsto |m|/|\tanh(m)|$, we obtain
\[
 \frac{1}{|\tanh(m)|} \le  \frac{1}{|\tanh(\delta)|}
\]
for $|m| \ge \delta$, and 
\begin{align*}
&\frac{\varepsilon}{|\tanh(m)|} + \frac{(1 - \varepsilon )|m|}{\delta |\tanh(m)|} 
\;\le\;  \frac{\varepsilon}{|\tanh(\varepsilon)|} + \frac{(1 - \varepsilon )|\delta|}{\delta |\tanh(\delta)|} \;\le\; \frac{1+\delta}{\tanh( \delta )}
\end{align*}
for $\varepsilon \le |m| \le \delta$. This yields
\begin{align}
\sup_{m \in \R}\, |\widehat{K}_0^{-1}(m)\, \widehat{\vartheta}(m)\, \widehat{P}_{\varepsilon, \infty}(m)| \leq \frac{1+\delta}{\tanh( \delta )}\,.
\end{align}
Furthermore, we have
\begin{align}
\sup_{k-m \in \R}\, |\widehat{K}_0^{-1}(k-m)\,  \chi_c(k-m)| \leq \frac{1}{\tanh(k_0-\delta )}\,.
\end{align}
The definitions of $\vartheta$ and $P_{\delta,\infty}$ directly imply
\begin{align} \label{theta-1-eps-1}
\sup_{k \in \R}\, |\widehat{\vartheta}^{-1}(k)| &= \varepsilon^{-1}\,,\\
\sup_{k \in \R}\, |\widehat{P}_{\delta,\infty}(k)\widehat{\vartheta}^{-1}(k)| &= 1\,.
\end{align}
Moreover, we have
\begin{align*}
|k\, \widehat{\vartheta}^{-1}(k)|= 
\begin{cases}
|k| & \quad {\rm for }\; |k| > \delta\,,\\[1mm]
\dfrac{|k|}{\varepsilon +(1- \varepsilon) \frac{|k|}{\delta}} & \quad{\rm for}\; |k| \le \delta \,.
\end{cases}
\end{align*}
Since 
\begin{align*}
\frac{|k|}{\varepsilon +(1- \varepsilon) \frac{|k|}{\delta}} \;=\; 
\frac{1}{\frac{\varepsilon}{|k|} + \frac{(1- \varepsilon)}{\delta}} \;\le\; 
\frac{1}{\frac{\varepsilon}{\delta} + \frac{(1- \varepsilon)}{\delta}}\;=\; \delta
\end{align*}
for $0 \neq |k| \le \delta$, we get
\begin{align}
\label{d/dx theta-1}
\sup_{k \in \R}\, |k\widehat{\vartheta}^{-1}(k)| &= \max\{\delta, |k|\}\,.
\end{align}
Now, using \eqref{NFourier}-\eqref{d/dx theta-1}, \eqref{kthk}, \eqref{RES3}, Young's inequality for convolutions,
\begin{align*}
\widehat{n}(k,k-m,m)=\widehat{n}(-k,-(k-m),-m) \in \R
\end{align*}
and the fact that $\psi_c$ is real-valued, we obtain the validity of all statements of a).

Let  $k_1>0$ be a constant such that $|k| \ge k_1$ and $|k-m- k_0 | \le \delta$ imply
$|m| \ge \delta$ and $\sign(k)=\sign(m)$.
Then, by using 
\begin{align*}
\tanh(k)=\sign(k) \Big(1-\frac{2}{1+e^{2 |k|}}\Big)
\end{align*}
we get 
\begin{align*}
\widehat{\vartheta}(k)\, \widehat{n}(k,k-m,m)&=-\frac{k\, \chi_c(k-m)}{\tanh(k-m)}(1+\mathcal{O}(e^{-2|k|}))
(1+\mathcal{O}(e^{-2|k-(k-m)|})) \\[2mm]
&= \frac{k\, \chi_c(k-m)}{\tanh(k-m)}(1+\mathcal{O}(e^{-|k|}))
\end{align*}
for $|k| \ge k_1$ provided that $k_1$ is chosen large enough. This yields statement b).

\eqref{N-wahl} follows by construction of $N$ due to \eqref{k0-id}.
\eqref{P_0 N} is a direct consequence of
\begin{align*}
\chi_{[-\delta, \delta]}(k) \chi_{[-\delta, \delta]}(m) \chi_c(k-m)=0\,.
\end{align*}
Finally, \eqref{partN} follows from a) and b) by  integration by parts.
\qed

\begin{lemma}
\label{int-kerne}
Fix $p \in \R$. Assume that $\kappa \in C(\R^3,\C)$, that $g \in C^2(\R,\C)$ has a compactly supported Fourier transform and that
$f \in H^{s}(\R,\C)$ for $s \geq 0$.\\[2mm]
{\bf a)} If $\kappa$ is Lipschitz continuous with respect to its second argument in some neighborhood
of $p$, then there exist $C_{g,\kappa,p} >0$, $\eps_0 >0$ such that
\begin{equation} 
\Big \| \int \big(\kappa (\cdot, \cdot-\ell,\ell) - \kappa (\cdot, p ,\ell)\big)\, \veps^{-1} \widehat{g}\Big(\frac{ \cdot - \ell - p}{\varepsilon}\Big) \widehat{f}(\ell)\, d\ell\, \Big \|_{L^{2}(s)}
 \le C_{g,\kappa,p}\, \varepsilon \| f \|_{H^{s}}
\end{equation}
for all $\eps \in (0, \eps_0)$.
\\[2mm]
{\bf b)} If $\kappa$ is globally Lipschitz continuous with respect to its third argument, then there exist
$D_{g,\kappa} >0$, $\eps_0 >0$ such that 
\begin{equation}  
\Big \| \int \big(\kappa (\cdot, \cdot-\ell,\ell)- \kappa(\cdot, \cdot - \ell ,\cdot-p)\big)\, \veps^{-1} \widehat{g}\Big(\frac{ \cdot - \ell - p}{\varepsilon} \Big) \widehat{f}(\ell)\, d\ell\, \Big \|_{L^{2}(s)}
\le D_{g,\kappa}\, \varepsilon \| f \|_{H^{s}} 
\end{equation}
for all $\eps \in (0, \eps_0)$.
\end{lemma}

\textbf{Proof.}
The Lemma is a special case of Lemma 3.5 in \cite{DSW12}.
\qed

\begin{lemma}
\label{Tlem}
The operator $\mathcal{T}$ has the following properties:\\[2mm]
{\bf a)}
Fix functions $g,h$ with $\widehat{g}_c:= \chi_{\,\supp(\psi_c)}\, \widehat{g} \in L^1(\mathbb{R}, \C)$ and $ \hat{h}_c:= \chi_{\,\supp(\psi_c)}\, \widehat{h} \in L^1(\mathbb{R}, \C)$.
Then $f \mapsto \mathcal{T}(g_c,h_c,f)$ defines a continuous linear map from $L^2(\R,\C)$ into $L^2(\R,\C)$, and there exists a constant $C>0$ such that for all $f \in L^2(\R,\C)$ we have
\begin{align}
\label{T-abschaetzung}
\Vert \mathcal{T}(g_c,h_c,f) \Vert_{L^2} \le C {\varepsilon}^{-1} \Vert \widehat{g}_c\Vert_{L^1} \Vert \widehat{h}_c\Vert_{L^1} \Vert  f\Vert_{L^2}\,.
\end{align}
{\bf b)}
For all $f \in H^2(\R,\C)$ we have 
\begin{align}
\label{T-wahl}
&-\k0 \mathcal{T}(\psi_c,\psi_c,f)+
\mathcal{T}(\k0\psi_c,\psi_c,f) +\mathcal{T}(\psi_c,\k0\psi_c,f) 
+\mathcal{T}(\psi_c,\psi_c,\k0 f) \\ \nonumber
&\qquad\quad = N(\psi_c, \vartheta^{-1} \partial_x(\psi \vartheta f)) + Y(\psi,f)
\end{align}
with
\begin{equation}
\label{T-rest}
\|Y(\psi,f)\|_{L^2} = \mathcal{O}(\Vert f \Vert_{H^2})
\end{equation}
for sufficiently small $\eps >0$.\\[2mm]
{\bf c)}
For all $f \in L^2(\R,\C)$ we have
\begin{align}
\label{P_0 T}
P_{\delta,\infty}T(\psi_c, \psi_c,f)=0\,.
\end{align}
\end{lemma}

\textbf{Proof.}
To show a), we use the triangle inequality, Young's inequality for convolutions, \eqref{kthk}, \eqref{theta-1-eps-1} and \eqref{delta} to get
\begin{align*}
\Vert T(g_c,h_c,f) \Vert_{L^2} &\lesssim \sum_{j= \pm 1} \|\widehat{t}_j\|_{L^{\infty}}  \Vert \widehat{g}_c\Vert_{L^1} \Vert  \widehat{h}_c\Vert_{L^1} \Vert  f\Vert_{L^2} \lesssim  \varepsilon^{-1} \Vert \widehat{g}_c\Vert_{L^1} \Vert  \widehat{h}_c\Vert_{L^1} \Vert  f\Vert_{L^2} \,.
\end{align*}

To prove b), we first show that
\begin{align}
\label{Ndom}
N(\psi_c, \vartheta^{-1} \partial_x(\psi \vartheta f)) = \sum_{j= \pm 1} P_{0, \delta} N(\psi_j, \vartheta^{-1} \partial_x(\psi_j \vartheta f)) + \mathcal{O}(\Vert f \Vert_{H^2})
\end{align}
such that it is sufficient to prove that the $L^{2}$-norm of 
\begin{align*}
\tilde{Y}:=& \sum_{j= \pm 1} \Big( -\k0 \mathcal{T}_j(\psi_j,\psi_j,f)+
\mathcal{T}_j(\k0\psi_j,\psi_j,f) + \mathcal{T}_j(\psi_j,\k0\psi_j,f)\\
&\qquad\quad +\mathcal{T}_j(\psi_j,\psi_j,\k0 f) - P_{0, \delta} N(\psi_j, \vartheta^{-1} \partial_x(\psi_j \vartheta f)) \Big)
\end{align*}
is of order $\mathcal{O}(\Vert f \Vert_{L^2})$, which we will obtain
by construction of $\mathcal{T}$ and because of Lemma \ref{int-kerne}.

To verify \eqref{Ndom}, we split $N$ into
\begin{align*}
N(\psi_c, \vartheta^{-1} \partial_x(\psi \vartheta f)) =& \sum_{j= \pm 1} P_{0, \delta} N(\psi_j, \vartheta^{-1} \partial_x(\psi_j \vartheta f)) +
\sum_{j= \pm 1} P_{0, \delta} N(\psi_j, \vartheta^{-1} \partial_x(\psi_{-j} \vartheta f))\\
&+ P_{\delta, \infty}N(\psi_c, \vartheta^{-1} \partial_x(\psi \vartheta f))
+ \eps P_{0, \delta}N(\psi_c, \vartheta^{-1} \partial_x(\psi_s \vartheta f))\,.
\end{align*}
Due to \eqref{N0}, \eqref{PN} and \eqref{d/dx theta-1}, the $L^{2}$-norm of the sum of the last two summands is of order $\mathcal{O}(\Vert f \Vert_{H^2})$.
Furthermore, in Fourier space, we have
\begin{align*}
&\;\widehat{P_{0, \delta}N}(\psi_j, \vartheta^{-1} \partial_x(\psi_{\ell} \vartheta f))(k) \\[2mm]
=&\;
\widehat{P}_{0, \delta}(k) \int_{\mathbb{R}}\int_{\mathbb{R}} {K}(k,k-m,m,n)
\widehat{\psi}_{j}(k-m) \widehat{P}_{\eps,\infty}(m) \widehat{\psi}_{\ell}(m-n) \widehat{f}(n)\, dn dm
\end{align*}
with
\begin{align*}
{K}(k,k-m,m,n)=-\frac{ikm\,\widehat{\vartheta}(n)}{\widehat{\vartheta}(k)
\tanh(k)  \tanh(k-m) \tanh(m) }\,,
\end{align*}
where $y \mapsto y/\tanh(y)$ is continued by $1$ for $y=0$.

For $\ell=-j$ we can apply Fubini's theorem, Young's inequality for convolutions and Lemma \ref{int-kerne} to obtain 
\begin{align*}
&\;\|\widehat{P_{0, \delta}N}(\psi_j, \vartheta^{-1} \partial_x(\psi_{-j} \vartheta f))\|_{L^{2}}\\[2mm]
=&\;
\Big \| \int_{\mathbb{R}} \int_{\mathbb{R}} \widehat{P}_{0, \delta}(\cdot)\, {K}(\cdot,jk_0,\cdot-jk_0,\cdot)
\widehat{\psi}_{j}(\cdot-m) \widehat{P}_{\eps,\infty}(m) \widehat{\psi}_{-j}(m-n)\widehat{f}(n)\, dn dm \,\Big \|_{L^{2}}\\[2mm]
&\; + \mathcal{O}(\Vert f \Vert_{L^2})\,.
\end{align*}
Since
\begin{align*}
{K}(k,jk_0,k-jk_0,k) =-\frac{ik(k-jk_0)\, \widehat{\vartheta}(k)}{\widehat{\vartheta}(k)
\tanh(k) \tanh(jk_0)\tanh(k-jk_0)}
\end{align*}
and the factor $\widehat{\vartheta}(k)$ in the denominator is canceled by the same factor in the numerator, ${K}(k,jk_0,k-jk_0,k)$ contains no factors which are of order $\mathcal{O}(\eps^{-1})$ such that
\begin{align*}
\sup_{k \in \R}\, |{K}(k,jk_0,k-jk_0,k)| = \mathcal{O}(1)\,.
\end{align*} 
Hence, by using \eqref{RES3} and Young's inequality for convolutions, we obtain 
\begin{align*}
\|\widehat{P_{0, \delta}N}(\psi_j, \vartheta^{-1} \partial_x(\psi_{-j} \vartheta f))\|_{L^{2}} =&\,
 \mathcal{O}(\Vert f \Vert_{L^2})
\end{align*}
such that we have verified
\eqref{Ndom}.

To estimate $ \Vert \tilde{Y} \Vert_{L^2}$, we use
\begin{align*}
&\;\widehat{\tilde{Y}}(k)\\[2mm]
=&\;
 \sum_{j= \pm 1} \widehat{P}_{0, \delta}(k) \int_{\mathbb{R}}\int_{\mathbb{R}} 
K_j(k,k-m,m-n,n)
\widehat{\psi}_j(k-m)\widehat{\psi}_j(m-n)\widehat{f}(n)\, dn dm 
\\[2mm]
&\;- \sum_{j= \pm 1} \widehat{P}_{0, \delta}(k) \int_{\mathbb{R}}\int_{\mathbb{R}} 
K(k,k-m,m,n)
\widehat{\psi}_j(k-m)\widehat{P}_{\eps,\infty}(m)\widehat{\psi}_j(m-n)\widehat{f}(n)\, dn dm 
\,,
\end{align*}
where
\begin{align*}
{K}_j(k,k-m,m-n,n) =&\; -\frac{ik(k-jk_0)\,\widehat{\vartheta}(k-2jk_0)}{\widehat{\vartheta}(k) 
\tanh(k) \tanh(jk_0) \tanh(k-jk_0)} \\[2mm] 
&\; \times \frac{\tanh(k)-\tanh(k-m)-\tanh(m-n)-\tanh(n)}{\tanh(k) -2 \tanh(jk_0)-\tanh(k-2jk_0)}
\end{align*}
and $K$ is as above. We can apply again Fubini's theorem, Young's inequality for convolutions and Lemma \ref{int-kerne} to obtain
\begin{align*}
&\,\widehat{\tilde{Y}}(k)\\[2mm]
=&\,
 \sum_{j= \pm 1} \int_{\mathbb{R}}\int_{\mathbb{R}} 
\widehat{P}_{0, \delta}(k)\, K_j(k,jk_0,jk_0,k\!-\!2jk_0)
\widehat{\psi}_j(k-m)\widehat{P}_{\eps,\infty}(m)\widehat{\psi}_j(m\!-\!n)\widehat{f}(n)\, dn dm 
\\[2mm]
&\,-\!\sum_{j= \pm 1} \int_{\mathbb{R}}\int_{\mathbb{R}} 
\widehat{P}_{0, \delta}(k)\, K(k,jk_0,k\!-\!jk_0,k\!-\!2jk_0)
\widehat{\psi}_j(k-m)\widehat{P}_{\eps,\infty}(m)\widehat{\psi}_j(m\!-\!n)\widehat{f}(n)\, dn dm 
\\[2mm]
&\,+ \mathcal{O}(\Vert f \Vert_{L^2})\,,
\end{align*}
where we used that due to the support of $\widehat{\psi}_j$ the first integrand vanishes for $|m| \leq \eps$. 
Since 
\begin{align*}
K_j(k,jk_0,jk_0,k-2jk_0) = K(k,jk_0,k-jk_0,k-2jk_0)\,,
\end{align*}
the two integral kernels, which are both of order $\mathcal{O}(\eps^{-1})$, cancel each other out such that we get
\begin{align*}
\Vert \tilde{Y} \Vert_{L^2} =  \mathcal{O}(\Vert f \Vert_{L^2})\,.
\end{align*}
Hence, we have proven b).

Finally, c) follows directly by the definition of $\mathcal{T}$.
\qed
\medskip

\begin{lemma}
\label{int-helfer}
Let $f \in H^{\ell}(\R,\R)$ and $g \in H^{m}(\R,\R)$ with $\ell,m \ge 0$. Then we have 
\begin{align}
\label{int und theta}
\int_{\R}\partial_x^{\ell} f\, \partial_x^{m} \vartheta g\, dx = \int_{\R}\partial_x^{\ell} f\,  \partial_x^{m} g\,dx+ \mathcal{O}(\Vert f\Vert_{L^2} \Vert g\Vert_{L^2})\,,
\end{align}
\begin{align}
\label{int und theta-1}
\int_{\R}\partial_x^{\ell} f\, \partial_x^{m+1} \vartheta^{-1} g\, dx = \int_{\R}\partial_x^{\ell} f\,  \partial_x^{m+1} g\, dx+\mathcal{O}(\Vert f\Vert_{L^2} \Vert g\Vert_{L^2})\,.
\end{align}
\end{lemma}

\textbf{Proof.}
Using the definition of $\vartheta$, we get
\begin{align*}
\int_{\R}\partial_x^{\ell} f\, \partial_x^{m} \vartheta g\, dx &= \int_{\R}\partial_x^{\ell} f\,  \partial_x^{m} g \,dx + (-1)^{\ell}\! \int_{\R}f\, \partial_x^{\ell+m} P_{0, \delta}(\vartheta-1) g \,dx \,,\\[2mm] 
\int_{\R}\partial_x^{\ell} f\, \partial_x^{m+1} \vartheta^{-1} g\,dx &= \int_{\R}\partial_x^{\ell} f\,  \partial_x^{m+1}g\,dx + (-1)^{\ell}\! \int_{\R}f\, \partial_x^{\ell+m+1} P_{0, \delta}(\vartheta^{-1}-1) g\, dx \,,
\end{align*}
which yields \eqref{int und theta} and, due to \eqref{d/dx theta-1}, also \eqref{int und theta-1}.
\qed
\medskip

\begin{lemma}
\label{checkRlem}
For sufficiently small $\eps > 0$ there exist constants $C, \check{C} >0$ such that 
\begin{equation}
\label{checkR gegen R}
\Vert \check{R} \Vert_{L^2} 
\leq C \Vert R \Vert_{H^{1}}\,,
\end{equation}
\begin{equation}
\label{R gegen checkR}
\Vert R \Vert_{L^2} \leq  \check{C} \Vert \check{R} \Vert_{L^2}\,. 
\end{equation}

\end{lemma}

\textbf{Proof.}
Estimate \eqref{checkR gegen R} is a direct consequence of the estimates \eqref{N0} and \eqref{T-abschaetzung}.

To prove \eqref{R gegen checkR} we introduce
$R_0:=P_{0, \delta}R$, $\check{R}_0:=P_{0, \delta}\check{R}$, $R_1:=P_{ \delta, \infty}R$, $\check{R}_1:=P_{ \delta, \infty}\check{R}$
and split $R$, $\check{R}$ into $R=R_0+R_1$ and $\check{R}=\check{R}_0+\check{R}_1$.
Because of \eqref{P_0 N} and \eqref{P_0 T}, $R_0$ satisfies 
\begin{equation}
R_0 + \varepsilon^{2} \mathcal{T}(\psi_c, \psi_c, R_0) = \check{R}_0 - \varepsilon P_{0,\delta} N(\psi_c, R_1)- \varepsilon^{2} \mathcal{T}(\psi_c, \psi_c, R_1)\,.
\end{equation}
Using \eqref{N0} and \eqref{T-abschaetzung} yields
\begin{align}
\label{R0est}
\Vert R_0 \Vert_{L^2}
& \lesssim \Vert \check{R}_0 \Vert_{L^2} + \Vert R_1 \Vert_{L^2}
\end{align}
for sufficiently small $\eps > 0$. Moreover, $R_1$ satisfies 
\begin{equation}
R_1 + \varepsilon P_{\delta,\infty} N(\psi_c, R_1)
= \check{R}_1 - \varepsilon P_{\delta,\infty} N(\psi_c, R_0)\,.
\end{equation}
Multiplying this equation with $R_1$, integrating and using $P_{\delta,\infty} N= P_{\delta,\infty} \vartheta N$ as well as \eqref{PN} yields
\begin{align*}
\Vert R_1 \Vert_{L^2}^2 + \varepsilon \int_{\mathbb{R}} R_1\,  \vartheta N(\psi_c, R_1)\, dx  
\lesssim  (\Vert \check{R}_1 \Vert_{L^2} +  \varepsilon \Vert R_0 \Vert_{L^2}) \Vert R_1 \Vert_{L^2}\,.
\end{align*}
Because of \eqref{partN} and \eqref{R0est}, we get
\begin{align}
\label{R1est}
\Vert R_1 \Vert_{L^2}
& \lesssim \Vert \check{R}_1 \Vert_{L^2} +  \varepsilon \Vert R_0 \Vert_{L^2}  \lesssim \Vert \check{R} \Vert_{L^2}
\end{align}
and
\begin{align}
\label{R0est2}
\Vert R_0 \Vert_{L^2}
& \lesssim \Vert \check{R} \Vert_{L^2} 
\end{align} 
for sufficiently small $\eps > 0$. Combining \eqref{R1est} and \eqref{R0est2} yields \eqref{R gegen checkR}.
\qed
\medskip

The assertions of Lemma \ref{Nlem} and Lemma \ref{checkRlem} imply

\begin{lemma}
For $\ell \geq 1$, we have
\begin{align}
\tilde{E}_{\ell} = \frac{1}{2} \Vert \partial_x^{\ell} R \Vert_{L^2}^2 +   \epsilon\, \mathcal{O}(\Vert R \Vert_{H^{\ell}}^2) \,.
\end{align}
\end{lemma}

\begin{corollary}
\label{en-Aequi}
$\sqrt{\tilde{\mathcal{E}}_{s}}$ is equivalent to $\Vert R \Vert_{H^{s}}$
for all $s \ge 1$ if $\eps > 0$ is sufficiently small.
\end{corollary}

Now, we are prepared to estimate the time derivative of $\tilde{\mathcal{E}}_{s}$ for any sufficiently regular solution of \eqref{d/dt R}. We obtain

\begin{lemma}
For sufficiently small $\eps > 0$, we have
\begin{equation}
\label{dt E_0}
\frac{d}{dt} \tilde{E}_0  \lesssim
 \varepsilon^2( \tilde{\mathcal{E}}_{2}
  + \varepsilon^{1/2} \tilde{\mathcal{E}}_{2}^{3/2} +1)\,.
\end{equation}
\end{lemma}

\textbf{Proof.}
Because of \eqref{d/dt R} and \eqref{ch R} we get
\begin{align*}
\frac{d}{dt} \tilde{E}_0 = \int_{\mathbb{R}} \overline{\check{R}}\, \partial_t\check{R}\, dx +\int_{\mathbb{R}} {\check{R}}\, \overline{\partial_t\check{R}}\, dx 
\end{align*}
with
\begin{align*}
\partial_t\check{R} =&\, \k0 \check{R} + \varepsilon^{-5/2} \vartheta^{-1}\res(\eps \psi)\\[2mm]
&\, - \eps \big( \vartheta^{-1} \partial_x(\psi_c \vartheta P_{\eps, \infty} R) 
+\k0 N(\psi_c, R)- N(\k0 \psi_c, R)- N(R, \k0 \psi_c)\big)
\\[2mm]
&\, + \eps \big(N( \partial_t \psi_c- \k0 \psi_c, R) +
N
(\psi_c,\varepsilon^{-5/2}\vartheta^{-1} \res(\eps \psi)) \big)
\\[2mm]
&\,- \eps \big(\vartheta^{-1} \partial_x({\psi}_c \vartheta P_{0,\eps} R) + \vartheta^{-1} \partial_x((\tilde{\psi}-\psi_c) \vartheta R) \big)
\\[2mm]
&\,- \eps^{2} 
N
(\psi_c,\vartheta^{-1} \partial_x({\psi} \vartheta R)) 
\\[2mm]
&\,- \eps^{2} 
\big( \k0 \mathcal{T}(\psi_c,\psi_c,R) 
- \mathcal{T}(\k0 \psi_c,\psi_c,R)-\mathcal{T}(\psi_c,\k0 \psi_c,R)- \mathcal{T}
(\psi_c,\psi_c,\k0R)\big)
\\[2mm]
&\,+ \eps^{2} 
\big(\mathcal{T}(\partial_t \psi_c- \k0 \psi_c, \psi_c,R) +\mathcal{T}(\psi_c,\partial_t \psi_c- \k0 \psi_c,R)\big) 
\\[2mm]
&\,+ \eps^{2} 
\mathcal{T}(\psi_c,\psi_c,\varepsilon^{-5/2}\vartheta^{-1} \res(\eps \psi))
\\[2mm]
&\,- \eps^{3}  \mathcal{T}
(\psi_c,\psi_c,\vartheta^{-1} \partial_x(\tilde{\psi} \vartheta R))
- \frac{1}{2}\eps^{7/2} 
N
(\psi_c,\vartheta^{-1} \partial_x(\vartheta R)^{2}) 
\,,
\end{align*}
where $\tilde{\psi}=\psi + \frac{1}{2}\veps^{3/2} \vartheta R$.

Exploiting the skew symmetry of $\k0 \,$ and 
the Cauchy-Schwarz inequality, we conclude
\begin{align*}
\frac{d}{dt} \tilde{E}_0 
\leq   2 \Vert \partial_t\check{R} -  \k0 \check{R} \Vert_{L^2} \Vert \check{R} \Vert_{L^2}\,.
\end{align*}
Using \eqref{ans-high}, the bounds \eqref{RES1}, \eqref{RES3} and \eqref{d/dt Psi}
for the approximation functions and the residual, the properties \eqref{N0} and  
\eqref{N-wahl} of the operator $N$, the properties \eqref{T-abschaetzung}-\eqref{T-rest} of the operator $\mathcal{T}$, the bounds \eqref{d/dx theta-1} and
\begin{equation} 
\label{P_0 Theta} 
 \|\vartheta P_{0,\epsilon} f \|_{L^{2}} \lesssim \epsilon \|f \|_{L^{2}} 
\end{equation} 
for $\vartheta$, the estimate \eqref{checkR gegen R} for $\check{R}$ as well as Corollary \ref{en-Aequi}, we get
\begin{align*}
\Vert \partial_t\check{R} -  \k0 \check{R} \Vert_{L^2}
\lesssim &\; \veps^{2}  ( \tilde{\mathcal{E}}_{2}^{1/2}
  + \varepsilon^{1/2} \tilde{\mathcal{E}}_{2} +1) \,,\\[2mm]
\Vert \check{R} \Vert_{L^2}
\lesssim &\; \tilde{\mathcal{E}}_{1}^{1/2}
\end{align*}
and therefore
\begin{equation*}
\frac{d}{dt} \tilde{E}_0  \lesssim
 \varepsilon^2( \tilde{\mathcal{E}}_{2}
  + \varepsilon^{1/2} \tilde{\mathcal{E}}_{2}^{3/2} +1)\,.
\end{equation*}
\qed
\medskip

\begin{lemma}
For $\ell\ge 1$, $\theta \ge \max \{2,\ell\} $ and sufficiently small $\eps > 0$,  we have
\begin{equation}
\label{dt E_l}
\frac{d}{dt} \tilde{E}_{\ell}  \lesssim
 \varepsilon^2 (\tilde{\mathcal{E}}_{\theta} +
  \varepsilon^{1/2} \tilde{\mathcal{E}}_{\theta}^{3/2}
  	+1)\,.
\end{equation}
\end{lemma}

\textbf{Proof.}
We compute
\begin{align*}
\frac{d}{dt} \tilde{E}_{\ell} \;=&\; \int_{\R} \partial^{\ell}_x R\, \partial_t\partial^{\ell}_x R \,dx\,
+   \veps \Big( \int_{\R} \partial_t\partial^{\ell}_x R\, \partial^{\ell}_x N(\psi_c,R)\,dx
\\[2mm] &\;
+ \int_{\R} \partial^{\ell}_x R\, \partial^{\ell}_x N(\psi_c,\partial_t R)\,dx\,
+ \int_{\R} \partial^{\ell}_x R\, \partial^{\ell}_x N(\partial_t \psi_c, R)\,dx \Big)\,.
\end{align*}
Using the error equation \eqref{d/dt R}, we get
\begin{align*}
\frac{d}{dt} \tilde{E}_{\ell} 
\,\,=&\,\, \int_{\R} \partial^{\ell}_x R \k0 \partial^{\ell}_x R \,dx 
\\[2mm]
&\,\,
+ \int_{\R} \partial^{\ell}_x R\, \veps^{-5/2} {\vartheta}^{-1} \partial^{\ell}_x \res(\veps\psi)\,dx\, 
\\[2mm] 
&\,\,
+\, \veps \Big(- \int_{\R} \partial^{\ell}_x R\, \vartheta^{-1} \partial_x^{\ell+1}(\psi_c \vartheta P_{\varepsilon, \infty} R)\,dx 
\\[2mm]
&\,\,
\qquad\;
+ \int_{\R} \k0 \partial^{\ell}_x R\, \partial^{\ell}_xN(\psi_c, R)\,dx
\\[2mm]
&\,\,
\qquad\;
+ \int_{\R} \partial^{\ell}_x R\, \partial^{\ell}_x N(\psi_c, \k0 R)\,dx
\\[2mm]
&\,\,
\qquad\;
+ \int_{\R} \partial^{\ell}_x R\, \partial^{\ell}_x N(\k0 \psi_c, R)\,dx
\\[2mm]
&\,\,
\qquad\;
+ \int_{\R} \partial^{\ell}_x R\, \partial^{\ell}_x N(\partial_t \psi_c -\k0 \psi_c, R)\,dx\\[2mm] 
&\,\,
\qquad\;
+ \int_{\R} \veps^{-5/2} {\vartheta}^{-1} \partial^{\ell}_x \res(\veps\psi)\, \partial^{\ell}_x N(\psi_c,R)\,dx
\\[2mm]&\,\, 
\qquad\;
+ \int_{\R} \partial^{\ell}_x R\, \partial^{\ell}_x N(\psi_c,\veps^{-5/2} {\vartheta}^{-1} \res(\veps\psi))\,dx\\[2mm]
&\,\,
\qquad\;
- \int_{\R} \partial^{\ell}_x R\, \vartheta^{-1} \partial_x^{\ell+1}(\psi_c \vartheta P_{0,\varepsilon} R)\,dx 
\\[2mm]&\,\, 
\qquad\;
- \int_{\R} \partial^{\ell}_x R\, \vartheta^{-1} \partial_x^{\ell+1}((\tilde{\psi}-\psi_c) \vartheta R)\,dx\, \Big) 
\end{align*}
\begin{align*}
 & -\, \veps^{2}  \Big(\, \int_{\R}  \vartheta^{-1} \partial_x^{\ell+1}(\tilde{\psi} \vartheta R)\, \partial^{\ell}_x N(\psi_c, R)\,dx\\[1mm]
&\,
\qquad\;\,\,
+\int_{\R} \partial^{\ell}_x R\, \partial^{\ell}_x N(\psi_c,\vartheta^{-1} \partial_x (\tilde{\psi} \vartheta R))\,dx \,\Big) \\[2mm]
\qquad\, =:&\, \sum_{j=1}^{13} I_j\,
\end{align*}
where $\tilde{\psi}=\psi + \frac{1}{2}\veps^{3/2} \vartheta R$.

Because of the skew symmetry of $\k0\,$, the integral $I_1$ equals zero. 
Since the operator $N$ satisfies \eqref{N-wahl}, we have $$I_3+I_4+I_5+I_6=0.$$ 

Furthermore, by integration by parts we get
\begin{align*}
I_8 =  (-1)^{\ell-1} \int_{\R} \veps^{-3/2} {\vartheta}^{-1} \partial^{2\ell-1}_x \res(\veps\psi)\, \partial_x N(\psi_c,R)\,dx \,.
\end{align*}
Hence, the bounds \eqref{RES1}, \eqref{d/dxN},  \eqref{d/dx theta-1} and the Cauchy-Schwarz inequality directly yield
\begin{align*}
I_8 = \veps^{4}\,\mathcal{O}(\tilde{\mathcal{E}}_{\theta}+ 1)
\,.
\end{align*}
Similarly, one can derive
\begin{align*}
I_2 + I_9 = \veps^{2}\,\mathcal{O}(\tilde{\mathcal{E}}_{\theta}+ 1)
\,.
\end{align*}

To control $I_7$ we use \eqref{int und theta-1} and obtain
\begin{align*}
I_7 =  \veps \int_{\R} \partial^{\ell}_x R\, \partial^{\ell}_x \vartheta N(\partial_t \psi_c -\k0 \psi_c, R)\,dx + \veps\, \mathcal{O}(\|R\|_{L^2} \|\vartheta N(\partial_t \psi_c -\k0 \psi_c, R)\|_{L^2} )
\end{align*} 
such that \eqref{N aufloesen}, \eqref{Q0}, \eqref{PI}, \eqref{d/dt Psi} and the Cauchy-Schwarz inequality yield
\begin{align*}
I_7 = \veps^{2}\,\mathcal{O}(\tilde{\mathcal{E}}_{\theta})
\,.
\end{align*}
Moreover, by using 
\eqref{RES3}, 
\eqref{d/dx theta-1}, \eqref{P_0 Theta}, the fact that  $\widehat{P_{0,\veps} f}$ has compact support and the Cauchy-Schwarz inequality we get 
\begin{align*}
I_{10} =& \veps^{2}\,\mathcal{O}(\tilde{\mathcal{E}}_{\theta})
\,.
\end{align*}

Next, we analyze $I_{12}+I_{13}$. Due to \eqref{int und theta-1},  we have
\begin{align*}
I_{12}+ I_{13} 
=& -\veps^{2}  \Big(\, \int_{\R}  \vartheta^{-1} \partial_x^{\ell+1}(\tilde{\psi} \vartheta R)\, \partial^{\ell}_x \vartheta N(\psi_c, R)\,dx\\[1mm]
&\,
\qquad\;\,\,
+\int_{\R} \partial^{\ell}_x R\, \partial^{\ell}_x \vartheta N(\psi_c,\vartheta^{-1} \partial_x (\tilde{\psi} \vartheta R))\,dx \,\Big) \\[2mm]
& +\, \veps^{2}\,\mathcal{O}(\tilde{\mathcal{E}}_{\theta}+ \veps^{3/2} \tilde{\mathcal{E}}^{3/2}_{\theta})\,. 
\end{align*}
To extract all terms with more than $\ell$ spatial derivatives falling on $R$, we apply Leibniz's rule and get
\begin{align*}
I_{12}+ I_{13} 
=&-\veps^{2}\, \Big(\,\int_{\R}  \vartheta^{-1} \partial_x^{\ell+1}(\tilde{\psi} \vartheta R)\, \vartheta N(\psi_c, \partial^{\ell}_x R)\,dx\\[1mm]
&\qquad\;\;\; + \ell \int_{\R}  \vartheta^{-1} \partial_x^{\ell+1}(\tilde{\psi} \vartheta R)\,  \vartheta N(\partial_x \psi_c, \partial^{\ell-1}_x R)\,dx\\[2mm]
&\qquad\;\;\; +\int_{\R} \partial^{\ell}_x R \,\vartheta N(\psi_c, \vartheta^{-1} \partial_x^{\ell+1}(\tilde{\psi} \vartheta R) )\,dx \\[2mm]
&\qquad\;\;\; +\ell \int_{\R} \partial^{\ell}_x R\, \vartheta N(\partial_x\psi_c, \vartheta^{-1} \partial_x^{\ell}(\tilde{\psi} \vartheta R))\,dx \,\Big)\\[2mm]
&+\, \veps^{2}\,\mathcal{O}(\tilde{\mathcal{E}}_{\theta}+ \veps^{3/2} \tilde{\mathcal{E}}^{3/2}_{\theta})\,. 
\end{align*}
Because of \eqref{partN}, we obtain
\begin{align*}
I_{12}+ I_{13}  =& -\veps^{2} \Big( \int_{\R}  \vartheta^{-1} \partial_x^{\ell+1}(\tilde{\psi} \vartheta R)\, {S}(\partial_x \psi_c, \partial^{\ell}_x R)\,dx\\[1mm]
&\qquad\;\;\;+2\ell \int_{\R}  \vartheta^{-1} \partial_x^{\ell+1}(\tilde{\psi} \vartheta R)\, \vartheta N(\partial_x \psi_c, \partial^{\ell-1}_x R)\,dx\, \Big)\\[2mm]
&+\, \veps^{2}\,\mathcal{O}(\tilde{\mathcal{E}}_{\theta}+ \veps^{3/2} \tilde{\mathcal{E}}^{3/2}_{\theta})\,. 
\end{align*}
Using \eqref{N aufloesen}, \eqref{int und theta} and \eqref{int und theta-1} yields
\begin{align*}
I_{12}+ I_{13} 
=& -(2\ell+1)\, \veps^{2} \int_{\R} \kk \partial_x \psi_c\, \tilde{\tilde{\psi}}\, \partial^{\ell}_x R\, \partial^{\ell+1}_x R\,dx \\[2mm]
&+\, \veps^{2}\,\mathcal{O}(\tilde{\mathcal{E}}_{\theta}+\veps^{3/2} \tilde{\mathcal{E}}^{3/2}_{\theta}) 
\end{align*}
where $\tilde{\tilde{\psi}}=\psi + \veps^{3/2} \vartheta R$.
Finally, with the help of \eqref{PI}, we arrive at
\begin{equation*}
I_{12}+ I_{13} 
= \veps^{2}\,\mathcal{O}(\tilde{\mathcal{E}}_{\theta}+\veps^{3/2} \tilde{\mathcal{E}}^{3/2}_{\theta})\,.
\end{equation*}

Using again \eqref{PI}, \eqref{int und theta} and \eqref{int und theta-1} as well as \eqref{RES3} yields
\begin{equation*}
I_{11} 
= \veps^{2}\,\mathcal{O}(\tilde{\mathcal{E}}_{\theta}+\veps^{1/2} \tilde{\mathcal{E}}^{3/2}_{\theta})\,.
\end{equation*}

Hence, we obtain 
\begin{equation*}
\frac{d}{dt} \tilde{{E}}_{\ell} \lesssim   \varepsilon^2 (\tilde{\mathcal{E}}_{\theta} +
  \varepsilon^{1/2} \tilde{\mathcal{E}}_{\theta}^{3/2}
  	+1)\,. 
\end{equation*}
\qed
\medskip

Now, combining the estimates \eqref{dt E_0} and \eqref{dt E_l}, we obtain 
\begin{equation}
\frac{d}{dt} \tilde{\mathcal{E}}_{s} \lesssim  \varepsilon^2 (\tilde{\mathcal{E}}_{s}
  + \varepsilon^{1/2} \tilde{\mathcal{E}}_{s}^{3/2}+1)
\end{equation}
for $s=s_A > 7$ and sufficiently small $\varepsilon>0$.
Consequently, Gronwall's inequality yields the $\mathcal{O}(1)$-boundedness of $\tilde{\mathcal{E}}_{s}$ for $t \in [0, T_0/{\eps}^2]$. 
Due to Corollary \ref{en-Aequi} and estimate \eqref{RES2}, Theorem \ref{Theorem 2 NLS-APP} follows.
\qed

\bibliographystyle{alpha}

\begin{thebibliography}{99}
\bibitem{AS81} M.J.~Ablowitz, H.~Segur, Solitons and the inverse scattering transform, in: SIAM Studies in Applied Mathematics, vol. 4. SIAM, 1981.

\bibitem{CDS15}
M.~Chirilus-Bruckner, W.-P.~D\"ull, G.~Schneider,
NLS approximation of time oscillatory long waves for equations with quasilinear quadratic terms,
{Math. Nachr.} {288} (2-3) (2015) 158-166.

\bibitem{CS11}
C.~Chong, G.~Schneider,
Numerical evidence for the validity of the NLS approximation in systems with a quasilinear quadratic nonlinearity,
{ZAMM Z. Angew. Math. Mech.} {93} (9) (2013) 688-696.

\bibitem{D12}
W.-P.~D{\"u}ll, 
Validity of the Korteweg-de Vries Approximation for the Two-Dimensional Water Wave Problem in the Arc Length Formulation,
{Comm. Pure Appl. Math.}  {65} (3) (2012) 381-429.

\bibitem{D16}
W.-P.~D{\"u}ll, 
Justification of the Nonlinear Schr\"odinger approximation for a quasilinear Klein-Gordon equation, 
{Comm. Math. Phys.} {355} (3) (2017) 1189-1207.

\bibitem{DS06}
W.-P.~D{\"u}ll, G.~Schneider, 
Justification of the Nonlinear Schr\"odinger equation 
for a resonant Boussinesq model,
{Indiana Univ. Math. J.} {55} (6) (2006) 1813-1834.

\bibitem{DSW12}
W.-P.~D{\"u}ll, G.~Schneider, C.E.~Wayne, 
Justification of the Nonlinear Schr\"odinger equation for the evolution of gravity driven 2D surface water waves in a canal of finite depth, 
{Arch. Rat. Mech. Anal.} {220} (2) (2016) 543-602.

\bibitem{HIT14}
J.K.~Hunter, M.~Ifrim, D.~Tataru, 
Two dimensional water waves in holomorphic coordinates,  Comm. Math. Phys. {346} (2) (2016) 483-552.

\bibitem{HITW13} J.K.~Hunter, M.~Ifrim, D.~Tataru, T.K.~Wong, Long Time Solutions for a Burgers-Hilbert Equation via a Modified Energy Method,
Proc. Amer. Math. Soc. {143} (8) (2015) 3407-3412. 

\bibitem{IT14}
M.~Ifrim, D.~Tataru,
The lifespan of small data solutions in two dimensional capillary water waves, Arch. Ration. Mech. Anal. 225 (3) (2017) 1279-1346.

\bibitem{IT16}
M.~Ifrim, D.~Tataru,
Two dimensional water waves in holomorphic coordinates II: Global solutions, Bull. Soc. Math. France {144} (2) (2016) 369-394.

\bibitem{Kal88}
L.A.~Kalyakin, 
Asymptotic decay of a one-dimensional wave packet in a nonlinear
  dispersive medium,
{Sb. Math. } {60} (1988) 457-483.

\bibitem{Kat75}
T.~Kato, The Cauchy problem for quasi-linear symmetric hyperbolic systems,
{Arch. Rat. Mech. Anal.} {58} (1975) 181-205.

\bibitem{MN13}
N.~Masmoudi, K.~Nakanishi, 
Multifrequency NLS scaling for a model equation of gravity-capillary
  waves,
{Commun. Pure Appl. Math.} {66} (8) (2013) 1202-1240.

\bibitem{Schn98JDE}
G.~Schneider, 
Approximation of the Korteweg-de Vries equation by the Nonlinear Schr\"odinger equation,
{J. Differential Equations} {147} (1998) 333-354. 

\bibitem{Schn05}
G.~Schneider, 
Justification and failure of the nonlinear Schr\"odinger equation
  in case of non-trivial quadratic resonances,
{J. Differential Equations} {216} (2005) 354-386.

\bibitem{Schn11OWbuch}
G.~Schneider, 
The role of the Nonlinear
Schr\"odinger equation in
nonlinear optics. In: 
Oberwolfach seminars 42:
Photonic Crystals: Mathematical Analysis and Numerical Approximation
by D\"orfler, W., Lechleiter, A., Plum, M., Schneider, G. and Wieners, C..
Birkh\"auser, 2011.

\bibitem{SSZ15} 
G.~Schneider, D.A.~Sunny, D.~Zimmermann, 
The NLS approximation makes wrong predictions  for the 
water wave problem  in case of small surface tension and spatially periodic boundary conditions,
{J. Dynam. Differential Equations} {27} (3) (2015) 1077-1099. 

\bibitem{SW10}
G.~Schneider, C.E.~Wayne, 
Justification of the NLS approximation for a quasilinear water
  wave model,
 J. Differential Equations {251} (2011) 238-269.  

\bibitem{SW00}
G.~Schneider, C.E.~Wayne, 
{ The long wave limit for the water wave problem.
I. the case of zero surface tension,}
Comm. Pure Appl. Math. {53} (12) (2000) 1475-1535. 

\bibitem{Schn11} 
G.~Schneider,  Justification of the NLS approximation for the KdV equation using
the Miura transformation, Advances in Mathematical Physics (2011) 854719.


\bibitem{Sh85}
J.~Shatah, 
Normal forms and quadratic nonlinear Klein-Gordon equations,
{Comm. Pure Appl. Math.} {38} (1985) 685-696.

\bibitem{T14}
N.~Totz, 
{A justification of the modulation approximation to the 3D full water wave problem,} 
Comm. Math. Phys. {335} (1) (2015), 369-443.

\bibitem{TW11} N.~Totz, S.~Wu, A rigorous justification of the 
modulation approximation to the 2D full water wave problem, 
{Comm. Math. Phys.} {310} (3) (2012) 817-883.


\end{thebibliography}

{\bf Address of the authors:}\\[2mm]
Universit\"at Stuttgart\\
Institut f\"ur Analysis, Dynamik und Modellierung\\
Pfaffenwaldring 57\\ 
70569 Stuttgart\\ 
Germany\\[2mm]
E-mail:\\
duell@mathematik.uni-stuttgart.de\\
max.hess@mathematik.uni-stuttgart.de

\end{document}